\documentclass[a4]{article}

\usepackage{amssymb,amsthm,amsmath}
\usepackage{amsfonts}
\usepackage{geometry}
\usepackage{graphicx}

\DeclareMathOperator{\vertex}{vert}

\def\Sym{\mathop{\rm Sym}\nolimits}
\def\acl{\mathop{\rm acl}\nolimits}
\def\Aut{\mathop{\rm Aut}\nolimits}
\def\Th{\mathop{\rm Th}\nolimits}
\def\tp{\mathop{\rm tp}\nolimits}
\def\Wr{\mathop{\rm Wr}\nolimits}

\def\Ff{\mathbb F} 

\def\AGL{\mathop{\rm AGL}\nolimits}
\def\PG{\mathop{\rm PG}\nolimits}
\def\PGL{\mathop{\rm PGL}\nolimits}

\def\calF{{\cal F}}
\def\calM{{\cal M}}
\def\calP{{\cal P}}
\def\calN{{\cal N}}
\def\calZ{{\cal Z}}

\newtheorem{defi}{Definition}[section]
\newtheorem{theorem}[defi]{Theorem}
\newtheorem{definition}[defi]{Definition}
\newtheorem{lemma}[defi]{Lemma}
\newtheorem{proposition}[defi]{Proposition}

\newtheorem{claim}{Claim}

\newtheorem{remark}[defi]{Remark}
\newtheorem{example}[defi]{Example}
\newtheorem{question}[defi]{Question}

\begin{document}

\title{Reducts of structures and maximal-closed permutation groups}

\author{Manuel Bodirsky\thanks{
The research of both authors was partially supported by EPSRC grant EP/H00677X/1, and the second author by EP/K010692/1. Manuel Bodirsky also received funding from the European Research Council under the European Community's Seventh Framework Programme (FP7/2007-2013 Grant Agreement no.\ 257039).} , \\ 
   Institut f\"{u}r Algebra\\TU Dresden\\01062 Dresden\\Germany\\
    Manuel.Bodirsky@tu-dresden.de,
\and Dugald Macpherson$^*$,\\
School of Mathematics\\
	University of Leeds\\
	 Leeds LS2 9JT, UK,\\
h.d.macpherson@leeds.ac.uk} 

\maketitle

\begin{abstract}
Answering a question of Junker and Ziegler, we construct a countable first order structure which is not $\omega$-categorical, but does not have any proper non-trivial reducts, in either of two senses (model-theoretic, and group-theoretic). We also construct a strongly minimal set which is not $\omega$-categorical but has no proper non-trivial reducts in the model-theoretic sense.
\end{abstract}

\section{Introduction}
\label{sect:intro}

For an $\omega$-categorical structure $\calM$, there is a clear notion of {\em reduct}
 $\calM'$ of $\calM$; namely, a structure whose domain  is equal to the domain $M$ of
 $\calM$, and such that any  subset of $M^k$ (for any $k>0$) which is $\emptyset$-definable
 in $\calM'$ is also $\emptyset$-definable in $\calM$. Two reducts are viewed as equal if 
they have the same $\emptyset$-definable sets. By the Ryll-Nardzewski Theorem, the
 notion 
can also be viewed group-theoretically: the automorphism group of a reduct is a 
closed subgroup of $\Sym(M)$ which contains $\Aut(\calM)$, two reducts are equal if and only if they 
have the same automorphism group, and every closed subgroup of $\Sym(M)$ containing 
$\Aut(\calM)$ corresponds to a reduct. (Closure here is with respect to the topology of pointwise convergence -- see the remarks after the statement of Theorem~\ref{reducts}.)

A countably infinite structure is {\em homogeneous} (in the sense of Fra\"iss\'e) if any isomorphism between finite substructures extends to an automorphism. By the Ryll-Nardzewski Theorem, any homogeneous structure over a finite relational language is $\omega$-categorical. 
For various such homogeneous structures, the reducts have been completely classified. See for example \cite{benn,bp1,cam1,thom1,thom2}.
Thomas has conjectured that if $\calM$ is homogeneous over a finite relational language then it has just finitely many reducts. Such questions have received additional motivation in recent work of the first author and collaborators on constraint satisfaction problems with infinite domains, where one aims to understand reducts up to primitive positive interdefinability (see, e.g.,~\cite{bp1}).

It appears that structures which are not $\omega$-categorical typically have infinitely many reducts. Junker and Ziegler \cite{jn}  asked whether this is always the case. In this paper, we give a negative answer to this question (see Theorem~\ref{reducts}), and construct a structure which is not $\omega$-categorical and has no proper non-trivial reducts. The example is a set equipped with a {\em $D$-relation} in the sense of \cite{an}. We also investigate the question for strongly minimal sets which are not $\omega$-categorical -- see Section 5 and in particular Theorem~\ref{dist-trans}.

Without the assumption of $\omega$-categoricity, the 
group-theoretic definition of `reduct' above does not in general coincide with the 
model-theoretic definition. We shall say that $\calM'$ is a {\em definable reduct} of 
$\calM$ if ${\calM}$ and $\calM'$ have the same domain $M$,  and every  subset of $M^k$ (for any $k$) which is $\emptyset$-definable in $\calM'$ is 
$\emptyset$-definable in $\calM$; two definable reducts are identified if they are inter-definable 
over $\emptyset$. A {\em group-reduct} of $\calM$ is a structure $\calM'$ with domain $M$ such that 
 $\Aut(\calM')$ is a closed subgroup of $\Sym(M)$ containing $\Aut(\calM)$, 
two group-reducts identified if their automorphism groups are equal. A definable reduct 
$\calM'$ of $\calM$ is {\em improper} if $\calM$ is $\emptyset$-definable in $\calM'$, and 
is {\em trivial} if $\calM'$ is $\emptyset$-definable in the pure set $M$, and we 
talk similarly of {\em improper} and {\em trivial} group-reducts. Clearly, any definable reduct is also a group-reduct, though a definable {\em proper} reduct may not be a proper group reduct. Also, if two definable reducts are different as group-reducts, then they are also distinct definable reducts.

We remark that it is also possible to mix the two notions, and consider group-reducts up to interdefinability. Considering reducts in this sense, the answer to the Junker-Ziegler question is positive -- see Proposition~\ref{mix} below. See also Proposition~\ref{posetemb} relating the two notions of reduct in the saturated case.

We now describe our example of Theorem~\ref{reducts}. If $D(x,y;z,w)$ is a quaternary relation on a set $M$, we say it is a {\em $D$-relation} if is satisfies the following axioms, taken from \cite{an}.
\begin{enumerate}
\item[(D1)] $\forall x \forall y \forall z \forall w \, \big( D(x,y;z,w) \to
( D(y,x;z,w) \wedge D(x,y;w,z) \wedge D(z,w;x,y)) \big)$
\item[(D2)] $\forall x \forall y \forall z \forall w \, \big( D(x,y;z,w) \to \neg D(x,z;y,w)\big)$
\item[(D3)] $\forall x \forall y \forall z \forall w \forall u \, \big(D(x,y;z,w) \to D(u,y;z,w) \vee D(x,y;z,u) \big)$
\item[(D4)] $\forall x\forall y \forall z \, \big ((x\neq z \wedge y \neq z) \to D(x,y;z,z) \big)$
\item[(D5)] $\forall x\forall y \forall z \, \big( x,y,z \mbox{~distinct}  \to \exists w \, (w\neq z \wedge D(x,y;z,w)) \big)$
\end{enumerate}

\begin{example}\label{mainex} \rm 
(1) {\em First description of example.} Let $M$ be the set of sequences of zeros and ones, indexed by ${\mathbb Z}$, which have {\em finite support}, that is, finitely many ones. Let $x,y,z\in M$,
with $x=(x_i)_{i\in {\mathbb Z}}$, $y=(y_i)_{i\in {\mathbb Z}}$ and $z=(z_i)_{i\in {\mathbb Z}}$.
Write $C(x;y,z)$ if there is $i\in {\mathbb Z}$ such that $x_i\neq y_i$ and such that for all $j\leq i$ we have $y_j=z_j$. Then $(M,C)$ is a {\em $C$-relation} in the sense of \cite{an}. Now define $D$ on 
$M$, putting $D(x,y;z,w)$ if and only if
one of the following holds: (a) $C(x;z,w)\wedge C(y;z,w)$, (b) $C(z;x,y) \wedge C(w;x,y)$. 

(2) {\em Second description of example.}
Let $(T,R)$ be the unique (graph-theoretic) tree of valency three, with vertex set $T$ and adjacency relation $R$. Let $M^+$ be the set of {\em ends} of $T$, that is, equivalence classes of one-way infinite paths (also called {\em rays}) of  $T$,
 where two paths  are equivalent if they have infinitely many common vertices. Let $x,y,z,w\in M^+$. Define $D(x,y;z,w)$ to hold if one of (a) $x=y \wedge x\neq z \wedge x\neq w$; (b) $z=w\wedge x\neq z \wedge y\neq z$; (c) $x,y,z,w$ are
 distinct and there are one-way paths $\hat{x}\in x, \hat{y}\in y,\hat{z}\in z,\hat{w}\in w$ such that $\hat{x} \cup \hat{y}$ and $\hat{z} \cup \hat{w}$ are disjoint two-way infinite paths (also called {\em lines}). Finally, let $M$ be a 
countable  subset of $M^+$ which is {\em dense}, in the sense that for any vertex $a\in T$ there are $x,y,z\in M$ and rays $\hat{x}\in x$, $\hat{y}\in y,\hat{z}\in z$ such that $\hat{x}\cup \hat{y}, \hat{x}\cup \hat{z}$ and $\hat{y} \cup \hat{z}$ are all lines through $a$.
\end{example}

For each of these descriptions, axioms (D1)--(D5) are easily verified -- see \cite[Theorem 23.5]{an} for description (1). It can be checked that the second description determines a unique structure up to isomorphism, and that the structures $(M,D)$ described in these two ways are isomorphic -- see Lemma~\ref{isom} below. We shall generally use the second description, as the symmetry is more visible. Until Section~\ref{sect:examples} of the paper, $\calM$ denotes the structure $(M,D)$. This structure is clearly not $\omega$-categorical. For example, as noted early in Section~\ref{sect:jordan} below, the tree $(T,R)$ is interpretable without parameters in $\calM$ and $\Aut(\calM)$ has infinitely many orbits on $T^2$ given by pairs at different distances, and infinitely many orbits on $M^4$. However, $\Aut(\calM)$ is 3-transitive on $M$ (see Lemma~\ref{isom}~$(iii)$). 

Our main theorem is the following, and gives a negative answer to Question 2 of \cite{jn}. 

\begin{theorem} \label{reducts}
\begin{enumerate}
\item[(1)] The structure $\calM$ has no proper non-trivial group-reducts.
\item[(2)] The structure $\calM$ 
 has no proper non-trivial  definable reducts.
\end{enumerate}
\end{theorem}

The proof of Theorem~\ref{reducts}~(1) is a straightforward application of classification results for primitive Jordan permutation groups from \cite{am}; see Definition~\ref{jordandef} below for the notion of {\em Jordan group}. The proof of (2) is more intricate: we show that in any non-trivial reduct the family of {\em cones} (see Section~\ref{sect:jordan}) is uniformly definable in the sense of Definition~\ref{unifdef} below, and hence that $\calM$ is definable, so the reduct is not proper.

\begin{definition} \rm \label{unifdef} If $\bar{x}=(x_1,\ldots,x_r)$ is a tuple, then we denote by $|\bar{x}|$ the length $r$ of $\bar{x}$. If $\calM$ is a structure, and $r\in {\mathbb N}$, we shall say that a family $\calF$ of subsets of $M^r$ is {\em uniformly definable} if there is $s\in {\mathbb N}$ and formulas 
$\phi(\bar{x},\bar{y})$ and $\psi(\bar{y})$, with $|\bar{x}|=r$ and $|\bar{y}|=s$, such that 
$$\calF:=\{\phi(M^r,\bar{a}): \bar{a}\in M^s \mbox{~and~} \calM\models \psi(\bar{a})\}.$$
\end{definition}

We give a group-theoretic interpretation of Theorem~\ref{reducts}(1).
If $M=\{a_n\: : \: n \in \omega\}$ is countably infinite, there is a natural complete metric $d$ on the symmetric group $\Sym(M)$: if $g,h\in \Sym(M)$ are distinct, put $d(g,h)=\frac{1}{n+1}$ if $n$ is least such that
$a_n^g\neq a_n^h$ or $a_n^{g^{-1}}\neq a_n^{h^{-1}}$, where following the conventions of this paper, $a_n^g$ denotes the image of $a_n$ under $g$. The group $\Sym(M)$ is a topological group with respect to the resulting topology, which is independent of the enumeration of $M$. A subgroup $G$ of $\Sym(M)$ is said to be {\em closed} if it is closed with respect to this topology. It is an easy exercise (see e.g. \cite[Section~2.4]{cam2}) to check that $G\leq \Sym(M)$ is closed if and only if there is a first order structure $\calM$ with domain $M$ such that $G=\Aut(\calM)$. 
We shall say that a closed proper subgroup $G$ of $\Sym(M)$ is {\em maximal-closed} if, whenever $G\leq H\leq \Sym(M)$ and $H$ is closed, we have $H=G$ or $H=\Sym(\calM)$. 

 Recall that a permutation group $G$ on a countably infinite set $M$ is {\em oligomorphic} if $G$ has finitely many orbits on $M^n$ for all $n>0$, or equivalently (by the Ryll-Nardzewski Theorem) if the topological closure of $G$ in $\Sym(M)$ is the automorphism group of an $\omega$-categorical structure with domain $M$. Thus Theorem~\ref{reducts} (1) provides an example of a maximal-closed proper subgroup of $\Sym(M)$ which is not oligomorphic;  maximal-closed subgroups previously known to us, such as those arising from the results cited earlier in \cite{cam1,thom1,thom2}, are oligomorphic. (We recently became aware of the results of \cite{bogomolov}. This yields other very different maximal-closed non-oligomorphic groups of the form $\PGL(n,{\mathbb Q})$ where $3\leq n\leq \aleph_0$, but we do not know if it yields other countable non-$\omega$-categorical structures without proper non-trivial definable reducts.) Possibilities for further non-oligomorphic examples are discussed below in Section~\ref{sect:examples}. We remark that, prior to the classification of finite simple groups (CFSG), proving that certain finite subgroups of the symmetric group $S_n$ were maximal was often technically  involved. In the infinite case, investigating maximal-closed subgroups is the natural analogue, and we have no version of CFSG available.

When looking for structures with no (or finitely many) reducts in either sense, it is natural 
to investigate classes of structures which are closed under  taking definable reducts and for which there is a good structure theory, or classes of closed permutation groups which are closed under taking closed supergroups (in the symmetric group) and have a good structure theory. The latter
consideration led us to $\calM$, since there is a structure theory for closed Jordan permutation groups. An obvious class of structures which is closed under taking definable reducts is the class of {\em strongly minimal sets}; that is, structures $\calM_1$ such that, for any $\calM_2\equiv \calM_1$, any definable subset of the domain of $\calM_2$ is finite or cofinite. 
We investigate this in Section~\ref{sect:examples}, and give an infinite family of examples (see Theorem~\ref{dist-trans}) of  strongly minimal structures which are 
not $\omega$-categorical but have no proper non-trivial definable reducts. We  have not been able to show that any such example has no proper non-trivial group-reducts  -- but note the recent work of Kaplan and Simon mentioned after Question~\ref{1,2}.

As further background, we briefly describe the first result of this kind for $\omega$-categorical structures, where the two notions of reduct coincide. Let $\calP$ be the homogeneous structure $({\mathbb Q},<)$. By \cite{cam1}, 
$\calP$  has three proper non-trivial reducts, namely $({\mathbb Q},B)$, $({\mathbb Q},K)$, and $({\mathbb Q},S)$. Here $B$ is the ternary 
{\em linear betweenness relation}
 on ${\mathbb Q}$, defined by putting
$B(x;y,z)$ whenever $y<x<z$ or $z<x<y$; the relation $K$ is the natural ternary {\em circular ordering} on ${\mathbb Q}$, where $K(x,y,z)$ holds if and only if $x<y<z$ or $y<z<x$ or $z<x<y$; and $S$ is the induced quaternary separation relation on ${\mathbb Q}$, where, for $x<y<z<w$ and
$\{s,t,u,v\}=\{x,y,z,w\}$,
we have $S(s,t;u,v)$ if and only if $\{s,t\}=\{x,z\}$  or
$\{s,t\}=\{y,w\}$. In particular, the $\omega$-categorical structure 
$({\mathbb Q},S)$ (itself a reduct of $({\mathbb Q},B)$ and $({\mathbb Q},K)$) has no proper non-trivial reducts; thus $\Aut({\mathbb Q},S)$ is maximal-closed in $\Sym({\mathbb Q})$.

Given the example in the last paragraph, a natural first attempt for a non-$\omega$-categorical structure with no proper non-trivial reduct (of either kind) would be $({\mathbb Z},S)$ where $S$ is the separation relation induced as above from the natural linear order on ${\mathbb Z}$. However, this has infinitely many distinct definable reducts. For example, it is easily seen that
the underlying graph on ${\mathbb Z}$ (with $u,v$ adjacent if and only if $|u-v|=1$) is $\emptyset$-definable in $({\mathbb Z},S)$, and has for every $n$ a reduct $\Gamma_n$, where $\Gamma_n$ is a graph on ${\mathbb Z}$ with $u,v$ adjacent if and only if $|u-v|=2^n$; these also yield infinitely many distinct group-reducts. See also Remark~\ref{morenotes}(3).

Finally, we mention the extensive model-theoretic literature  on reducts of fields, typically on structures which are definable reducts of a field $F$ in which the additive structure of the field (as  a module over itself) is $\emptyset$-definable. See for example \cite{marker} and 
\cite{peterzil}, which are related to Zilber Trichotomy phenomena.

Section~\ref{sect:jordan} below contains more background on  $D$-relations and the underlying example, and 
parts (1) and (2) of Theorem~\ref{reducts} are proved in Sections~\ref{sect:group-reducts} and \ref{sect:def-reducts}, respectively. Section~\ref{sect:examples} contains further discussion, mainly around reducts of strongly minimal sets, and some open questions.
We believe that the questions considered in this paper  have potential interest from both a model-theoretic and a permutation group-theoretic viewpoint, and have aimed to give sufficient background for readers from either perspective.

For space reasons the paper does not contain many diagrams, but the arguments in Sections 2--4 are highly pictorial, and we urge the reader to draw diagrams.

We thank the referee for a careful report which we believe has greatly benefitted the presentation of the paper. We also thank I. Kaplan for drawing our attention to \cite{bogomolov}.

 \section{Treelike structures and Jordan permutation groups}
 \label{sect:jordan}

We first discuss the structure $\calM$ rather more fully, and introduce some notation. 
Consider the tree $(T,R)$ from the second description of $\calM$, let $M^+$ be the set of all ends of $T$, and $M$ a countable dense subset of $M^+$, and $\calM=(M,D)$ (at this stage, we do not assume that this determines $\calM$ up to isomorphism -- see Lemma~\ref{isom} below).   We refer to elements of $T$ as {\em vertices} (of $T$, or of $M$). Given vertices $u,v$ of $T$ (or of any graph in the given context) we write $d(u,v)$ for the graph distance from $u$ to $v$. For each end $x\in M$ and $a\in T$ there is a unique ray $x_a\in x$ starting at $a$. We shall write $S(a,x_a)$, or just $S(a,x)$,  for the set of vertices on the tree $(T,R)$ lying on the ray $x_a$ (including $a$). There is an equivalence relation $E_a$ on $M$, where
 $E_axy$ holds if and only if the paths $x_a$ and $y_a$ have a common edge of $T$. The $E_a$-classes are called {\em cones} at $a$. Observe that as $(T,R)$ is trivalent,  there are three cones at $a$, that the union of two cones at $a$ is a cone at a neighbour of $a$, and (hence) that the complement in $M$ of a cone is a cone. For example, in Fig. 1, the union of the two leftmost cones at $a$ is a cone at $b$. If $U$ is a cone at $a$ we write $a=\vertex(U)$. 

\begin{figure}[h]
\begin{center}
 \includegraphics[scale=.5]{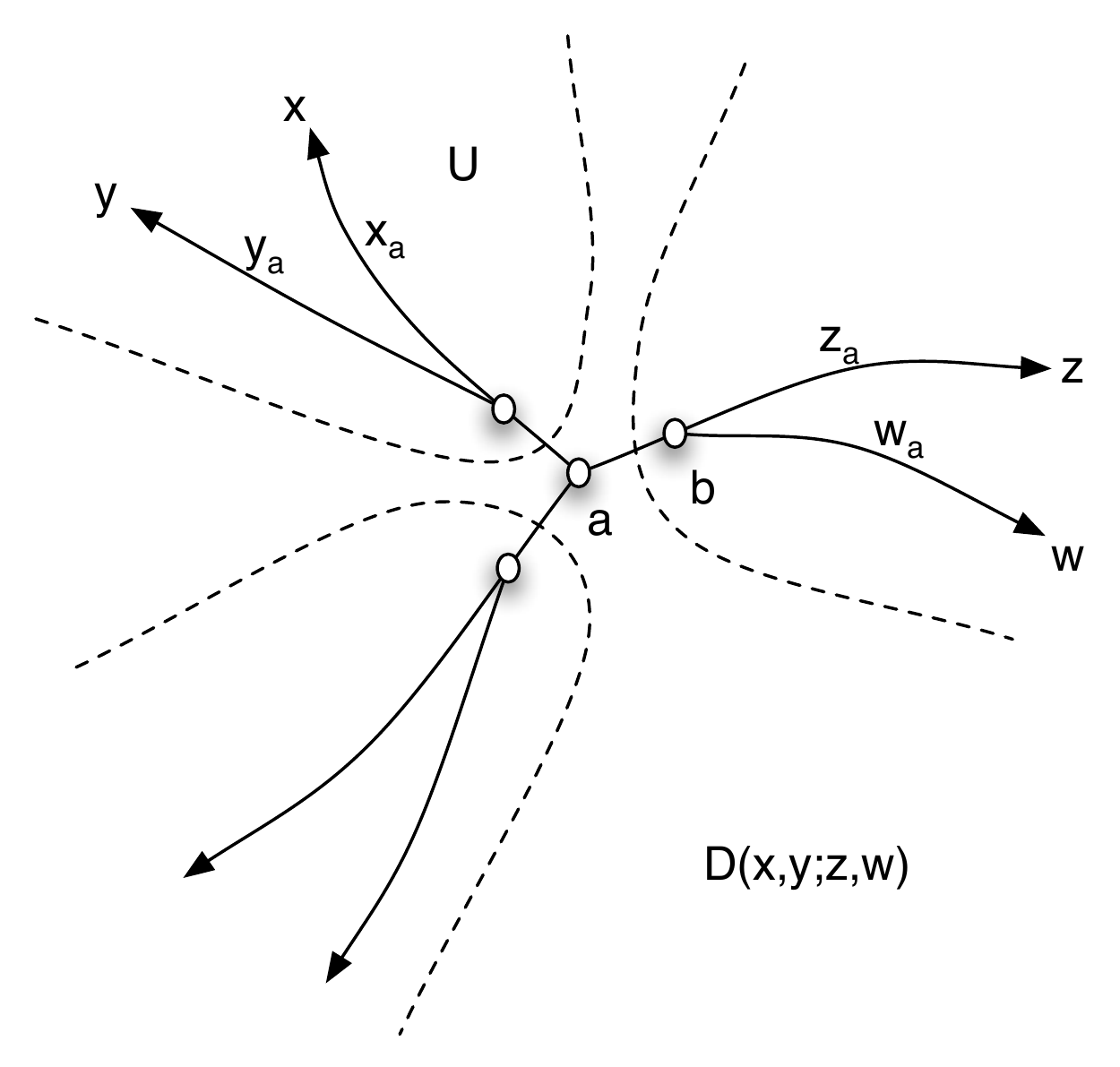}
\end{center}
\caption{Illustration for $\mathcal M$, $M^+$, and cones $U$ at $a$.}
\label{fig:xa}
\end{figure}

We may view a cone as a subset of $M$ or of $M^+$.
Slightly abusing notation, we sometimes view a cone $U$ at $a$ as a subset of $T$, namely as the union  $\bigcup_{x\in U}(S(a,x)\setminus\{a\})$. In particular, we may write $v\in U$ where $v$ lies in $T$ rather than $M$.

We note the following:

(i) any distinct $x,y,z\in M^+$ determine  a unique $a\in T$ such that $x,y,z$ are pairwise $E_a$-inequivalent;

(ii) any distinct $x,y\in M^+$ determine a unique line $l(x,y)$ (a connected 2-way-infinite path in the tree, viewed as a set of vertices); indeed, there exist $\hat{x} \in x, \hat{y} \in y$ such that $\hat{x} \cup \hat{y}$ is a connected 2-way infinite path in the tree, and $\hat{x} \cup \hat{y}$ is unique (it only depends on $x$ and $y$, not on the choice of the pair $\hat{x}, \hat{y}$). 

(iii) if $a\in T$ and $x,y\in M^+$ are distinct then $xE_ay$ if and only if $a\not\in l(x,y)$;

(iv) if $x,y,z\in M^+$ are pairwise $E_a$-inequivalent and $w\neq z$ then $wE_a z$ if and only if $l(x,y)\cap l(w,z)=\emptyset$;

(v) if $U_1,U_2$ are subcones of a cone $V$ then they are either disjoint or one contains the other.

\begin{figure}
\begin{center}
\includegraphics[scale=.5]{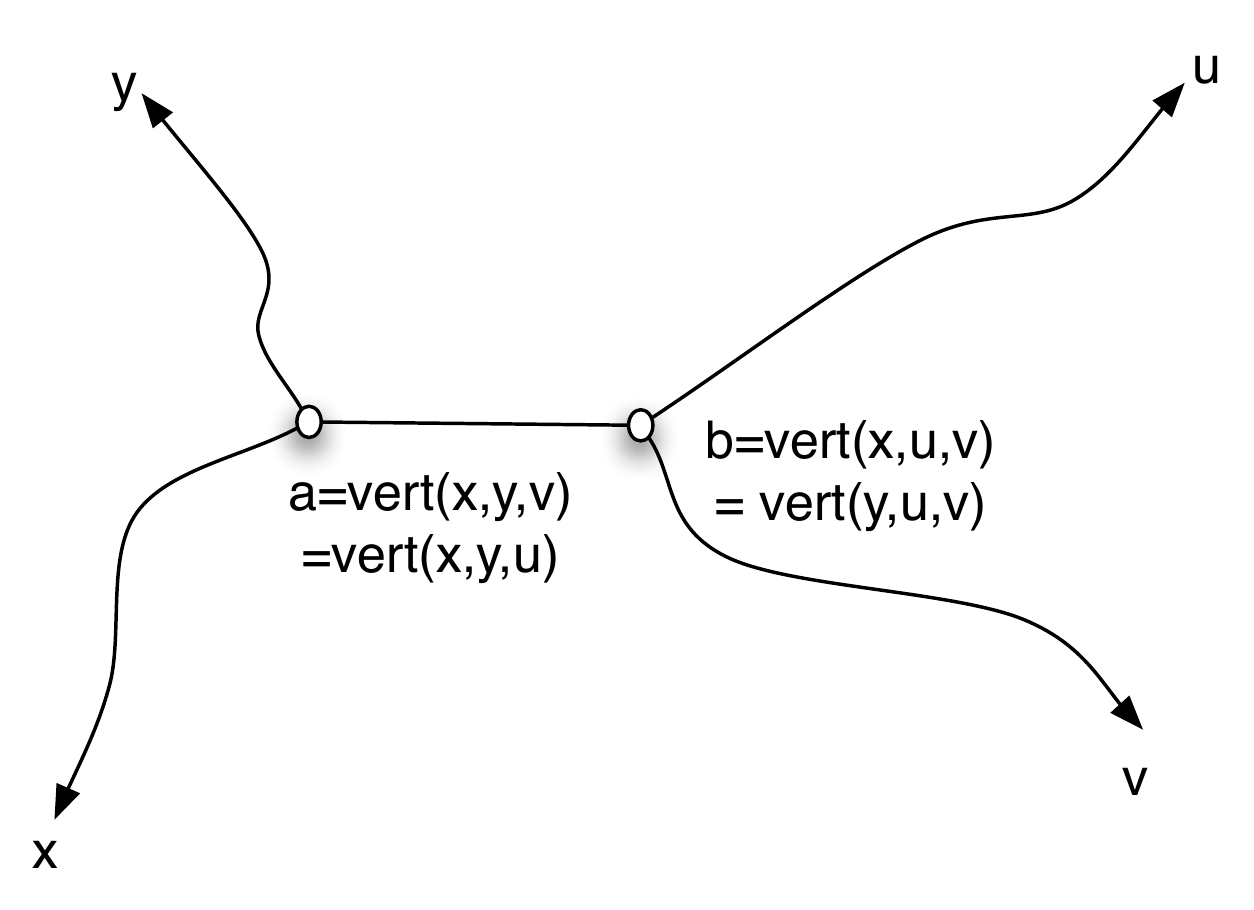}
\end{center}
\caption{Illustration for the interpretation of $(T,R)$ in $\calM$.}
\label{fig:vert}
\end{figure}

\begin{remark} \rm \label{cones}
\begin{enumerate}
\item  The set of all cones is uniformly definable in $\calM$. For if $x_a,y_a,z_a\in M$ lie in distinct cones at $a$, then the cone at $a$ containing $z_a$ is the set
 $\{w\in M: \calM\models D(x_a,y_a; z_a,w)\}$. Thus, if $\phi(w,xyz)$ is the formula
 $D(x,y;z,w)$ and $\psi(xyz)$ is the formula $x\neq y\wedge x\neq z\wedge y\neq z$, then the pair of formulas $\phi,\psi$ is a uniform definition of the family of cones, in the sense of Definition~\ref{unifdef}. 

The set of cones forms a basis of clopen sets for a totally disconnected topology on $M^+$ (and also on $M$). We say that a subset $X$ of $M^+$ is {\em dense} in $M^+$ if it is dense in this topology, that is, meets every cone. This coincides with the notion of density in Example~\ref{mainex}~$(2)$. 
 
 In \cite[Section~24]{an} cones are called {\em sectors}. These are defined more generally in \cite{an} for arbitrary $D$-relations -- see the discussion at the end of this section. In the proof of Theorem~\ref{reducts}~(1) we briefly refer to cones (i.e. sectors) for an arbitrary $D$-relation, not necessarily arising from $(T,R)$. 
 
\item The tree $(T,R)$ is interpretable in $\calM$. Indeed, 
as the set of cones is uniformly definable in $\calM$, we may define the 
set $S$ of all 3-sets
$\{U_1,U_2,U_3\}$ of cones such that the $U_i$ are pairwise disjoint and $U_1 \cup U_2 \cup U_3=M$. The set $S$ is then identified with $T$, identifying each vertex $a$ with the set of three cones at $a$. The vertex $\{U_1,U_2,U_3\}$ is adjacent to
the vertex $\{U_1',U_2',U_3'\}$ if there is  $U\in \{U_1,U_2,U_3\}$ and $U'\in \{U_1',U_2',U_3'\}$ such that $U'$ is the disjoint union of the two cones in $\{U_1,U_2,U_3\}\setminus U$, and $U$ is the disjoint union of the two cones in $\{U_1',U_2',U_3'\}\setminus U'$.
See also Sections 24--26 of \cite{an}.   

We may also view $T$ directly as a quotient of $M^3$. Let $Y$ be the set of all 3-element subsets of $M$. For $\{x,y,z\}, \{x',y',z'\}\in Y$, write
$E(\{x,y,z\},\{x',y',z'\})$ if there is $\{U_1,U_2,U_3\}\in S$ such that each $U_i$ contains exactly one of $x,y,z$ and exactly one of $x',y',z'$. 
Since for $\{x,y,z\}\in Y$ there is a unique $\{U_1,U_2,U_3\}\in S$ such that each $U_i$ contains exactly one of $x,y,z$, $E$ is an equivalence relation on $Y$.
Informally, $E(\{x,y,z\},\{x',y',z'\})$ says that the ends  $x,y,z$ and $x',y',z'$ `meet' at the same vertex. We may then identify $T$ with $Y/E$. If $x,y,z\in M$ all lie in distinct cones at the vertex $a$, we write $a=\vertex(x,y,z)$. We may define the graph relation $R$ on $T$ as follows: for $a,b\in Y/E$ we define $R(a,b)$ to hold if $a\neq b$ and there are distinct $x,y,u,v\in M$ with $D(x,y;u,v)$
such that $a=\vertex(x,y,u)=\vertex(x,y,v)$, $b=\vertex(x,u,v)=\vertex(y,u,v)$, and there is no $c\in Y/E\setminus\{a,b\}$  and $t\in M$ such that $c=\vertex(t,w,z)$ for all $w\in \{u,v\}$ and $z \in \{x,y\}$. See Figure~\ref{fig:vert}.

\item As noted before Theorem~\ref{reducts}, $\calM$ is not $\omega$-categorical, since it interprets $(T,R)$ which has pairs of vertices at arbitrary distance so is clearly not $\omega$-categorical. 
\end{enumerate}
\end{remark}

 Let $A\subset M$ be finite. We say that a vertex $v\in T$ is an {\em $A$-centre} if there are distinct $a,b,c\in A$ such that $v=\{a,b,c\}/E=\vertex(a,b,c)$ (in the notation of (2) above). A vertex of $T$ is an {\em $A$-vertex} if it is an $A$-centre or lies on  a path in $(T,R)$ between two $A$-centres. If $\bar{a}=(a_1,\ldots,a_n)$, then an {\em $\bar{a}$-centre}  or  {\em $\bar{a}$-vertex} is just an $\{a_1,\ldots,a_n\}$-centre or $\{a_1,\ldots,a_n\}$-vertex respectively. Any $A$-vertex is definable over $A$ in $\calM$ (as an element of $\calM^{{\rm eq}}$, formally). 

\begin{lemma} \label{isom}
Let $M_1$ and $M_2$ be countable dense subsets of $M^+$ in Example~\ref{mainex}~(2). Then
\begin{enumerate}

\item[(i)] For $i=1,2$ let $A_i$ be finite non-empty subsets of $M^+$, and let $T_i$ be the finite subtree  of $T$ consisting of all $A_i$-vertices 
(so if $|A_i|\leq 2$ then $T_i$ is empty). Let $f:(A_1,D)\to (A_2,D)$ be an isomorphism which induces an isomorphism $T_1 \to T_2$. Suppose $a_1\in M_1\setminus A_1$. Then there is $a_2\in M_2\setminus A_2$ such that $f$ extends to an isomorphism $(A_1\cup\{a_1\},D)\to (A_2\cup\{a_2\},D)$ which induces an isomorphism between the tree of $(A_1\cup\{a_1\})$-vertices and the tree of $(A_2\cup\{a_2\})$-vertices.
\item[(ii)] $\calM_1:=(M_1,D)$ and $\calM_2:=(M_2,D)$ are isomorphic.
\item[(iii)] If $\calM_1$ is as in $(ii)$, then $\Aut(\calM_1)$ is 3-transitive on $M_1$.
\item[(iv)] The structures $(M,D)$ described in Example~\ref{mainex} (1) and (2) are isomorphic.
\item[(v)] The group $\Aut(\calM)$ is transitive on the set of cones of $\calM$.
\end{enumerate}
\end{lemma}
\begin{proof}
$(i)$ We build  an isomorphism $f \colon \calM_1\to \calM_2$ by a back-and-forth argument.
At the $n^{{\rm th}}$ stage, we have finite $A_1\subset M_1$ and $A_2\subset M_2$, and an isomorphism $f_n \colon (A_1,D)\to (A_2,D)$. The sets $A_1$ and $A_2$ determine finite subtrees $T_1$ and $T_2$ of $T$, where $T_1$ consists of all $A_1$-vertices and $T_2$ of all $A_2$-vertices, 
and we suppose also that $f_n$ induces an isomorphism $f_n^T$ between $(T_1,R)$ and $(T_2,R)$. Now  let $a_1\in M_1\setminus A_1$, and $A_1':=A_1\cup\{a_1\}$. 

Assume first that $|A_1|>2$. Let $T_1'$ be the subtree of $T$ consisting of all $A_1'$-vertices. Note that 
 there is exactly one $A_1'$-centre $v_1=\vertex(a_1,b,c)$ (with $b,c\in A_1$) which is not an $A_1$-centre, 
and the vertices of $T_1'$ consist of those of $T_1$ together with (possibly) those lying between $v_1$ and the nearest $A_1$-centre. Easily there is a subtree $T_2'$ of $T$ such that $f_n^T$ extends to an isomorphism $T_1'\to T_2'$. Let $v_2:=f_n^T(v_1)$. As $\calM_2$ is dense, there is $a_2\in M_2$ such that $v_2=\vertex(a_2, f_n(b), f_n(c))$. Then extend $f_n$ to $f_{n+1}$ by putting $f_{n+1}(a_1)=a_2$.

If $|A_1|\leq 2$ then $T_1$ is empty. There are now no constraints over the choice of $a_2$, though if $A_1=\{b_1,c_1\}$ and $A_2=\{b_2,c_2\}$ then the resulting map $A_1\cup\{a_1\}\to A_2\cup\{a_2\}$ will induce a map taking $\vertex(a_1,b_1,c_1)$ to $\vertex(a_2,b_2,c_2)$.

$(ii), (iii)$ These are both immediate from (i).

$(iv)$ Let $M$ be the set of finite support sequences of zeros and ones indexed by ${\mathbb Z}$, with $D$ defined as in Examples~\ref{mainex}~(1). Define $T$ to be the collection of all subsequences of elements of $M$ indexed by a set of the form $(-\infty,n)$ for 
some $n\in {\mathbb  Z}$, and if $\sigma,\tau\in T$, define $R(\sigma,\tau)$ to hold if there is $n$ such that $\sigma$ is indexed by $(-\infty,n)$, $\tau$ by $(-\infty,n+1)$, and $\tau$ extends $\sigma$ (or vice versa).  It is easily seen that $(T,R)$ is a regular tree of valency 3, and that $M$ may be identified with a dense set of ends of $T$. An element $a\in M$ contains the ray consisting of the bounded  restrictions of $a$ (subsequences indexed by some $(-\infty,n)$), and consists of all rays which agree with this ray on some final segment.   Now apply $(ii)$. 

$(v)$ If $U_1,U_2$ are cones at vertices $a_1,a_2$ respectively, pick $x_i,y_i,z_i\in M$ such that $a_i=\vertex(x_i,y_i,z_i)$ and $z_i\in U_i$ (for each $i=1,2$). By (iii) there is $g\in \Aut(\calM)$ with $(x_1,y_1,z_1)^g=(x_2,y_2,z_2)$. We then have $U_1^g=U_2$.
\end{proof}

Consequences of Lemma~\ref{isom} (i) may be used without explicit mention. For example, suppose that $u,v,x,y\in M$ are distinct, as are $u',v',x',y'$, with $D(u,v;x,y)\wedge D(u',v';x',y')$. Suppose also that $a:=\vertex(u,v,x)$, $a':=\vertex(u',v',x')$, $b:=\vertex(u,x,y)$, $b':=\vertex(u',x',y')$, and that in the graph $T$ there is an $ab$-path $a=a_0,a_1,\ldots,a_n=b$ and an $a'b'$-path $a'=a_0',a_1',\ldots,a_n'=b'$. Suppose $a_i=\vertex(u,w_i,x)$ and $a_i'=\vertex(u',w_i',x')$, where $w_i,w_i'\in M$, for $i=1,\ldots,n-1$. Then there is $g\in \Aut(\calM)$ with $(u,v,w_1,\ldots,w_{n-1},x,y)^g=(u',v',w_1',\ldots,w_{n-1}',x',y')$.

We shall heavily use properties of $\Aut(\calM)$, in particular, the fact that it is a {\em  Jordan permutation group}.
If $G$ is a permutation group on $X$ (so denoted $(G,X)$), we say that $G$ is {\em $k$-transitive on $X$} if, for any distinct $x_1,\ldots,x_k\in X$ and distinct $y_1,\ldots,y_k\in X$ there is $g\in G$ with $x_i^g=y_i$ for $i=1,\ldots,k$; we say $G$ is {\em highly transitive} on $X$ if $G$ is $k$-transitive on $X$ for each positive integer $k$. The permutation group $(G,X)$ is
{\em  primitive} if there is no proper non-trivial $G$-invariant equivalence relation on $X$.
If $A\subset X$, we write $G_{(A)}$ for the {\em pointwise stabiliser} in $G$ of $A$, namely the group
$\{g\in G: g|_A= {\rm id}_A\}$. We put $G_{\{A\}}=\{g\in G: A^g=A\}$, the {\em setwise stabiliser} of $A$ in $G$.

\begin{definition} \label{jordandef} \rm Let $G$ be a transitive permutation group on a set $X$. A subset $A$ 
of $X$ is a {\em Jordan set} if $|A|>1$ and $G_{(X\setminus A)}$ is transitive on $A$. We say $A$ is a {\em proper} Jordan set if $A\neq X$, and, if $|X\setminus A|=n \in {\mathbb N}$, then $(G,X)$ is not $(n+1)$-transitive. A {\em Jordan group} is a transitive permutation group with a proper Jordan set.
\end{definition}

 \begin{remark} \rm \label{notes}
\begin{enumerate}
\item Every cone of $\calM$ is a Jordan set for $G:=\Aut(\calM)$. The main point to observe here is that (by Lemma~\ref{isom}(iv)) if $U$ is a cone, we may identify it with the set of sequences of zeros and ones indexed by ${\mathbb N}$ of finite support, carrying the ternary relation $C$, where $C(x;y,z)$ holds  if and only if either $y=z\neq x$, or the first index where $x$ and $y$ differ is before the first where $y$ and $z$ differ (see the  first description in Example~\ref{mainex} above). Any permutation of $U$ which preserves this $C$-relation extends to an element of $G_{(M\setminus U)}$. We may view $U$ as an elementary abelian 2-group (under pointwise addition), that is, an abelian group of exponent 2. The action of $U$ on itself by addition preserves the $C$-relation. Thus $U$ acts transitively by addition as a group of automorphisms of $(M,C)$, so $G_{(M\setminus U)}$ is transitive on $U$. It is also easily checked that if $x,y\in U$ are distinct then there is $g\in G_{(M\setminus U)}$ with $(x,y)^g=(y,x)$ -- let $g$ act on $U$ as addition by $x+y$. 
\item Since there are just three cones at each vertex, and since any cone is a Jordan set, we have the following.
 \begin{enumerate}
\item[(i)] Any cone at an $A$-vertex is $A$-definable in $\calM$. Indeed, suppose that $U$ is the cone at the $A$-vertex $u$. Then there are distinct $a,b\in A$ such that $u$ lies on the line $l(a,b)$. By
Remark~\ref{cones}(3) and the paragraph after the remark, the vertex $u$ and the set $W:=\{x:u=\vertex(a,b,x)\}$ are $A$-definable. The three cones at $u$ are $W$, the cone at $u$ containing $a$ which is defined as $D(a,M;b,x)$ for any $x\in W$, and the cone at $u$ containing $b$, defined by $D(b,M,a,x)$ for any $x\in W$. These definitions are over $A$.

\item[(ii)] If $v$ is an $A$-vertex and $U$ is a cone at $v$ not containing any element of $A$, then $U$ is an orbit of $\Aut(\calM)_{(A)}$.
\item[(iii)] If $v$ is an $A$-vertex, $a \in A$, and $u$ is a vertex of $S(v,a)$, then induction on $d(u,v)$ shows that $u$ is $A$-definable (as an element of $M^{{\rm eq}}$), as is each of the three cones at $u$.
\item[(iv)] If $A\subset M$ is finite with $|A|\geq 3$, then the orbits of $\Aut(\calM)_{(A)}$ which are disjoint from $A$ are of two kinds: (a) any cone which contains no member of $A$ and is at an $A$-vertex $w$ which is not an $A$-centre; and (b), for any $A$-centre $u$, cone $U$ at $u$, and $a\in A$ such that $U\cap A=\{a\}$, and any $v\neq u$ on $S(u,a)$ with $v$ not an $A$-centre, the unique cone at $v$ which is disjoint from $A$. See Figure~\ref{fig:orbits}. Note that in (b),  by (iii) such a vertex $v$ is $A$-definable, and by (1) such a cone is contained in an $\Aut(\calM)_{(A)}$-orbit so  {\em is} an orbit, as the other cones at $v$ meet $A$.  
\begin{figure}[h]
\begin{center}
 \includegraphics[scale=.5]{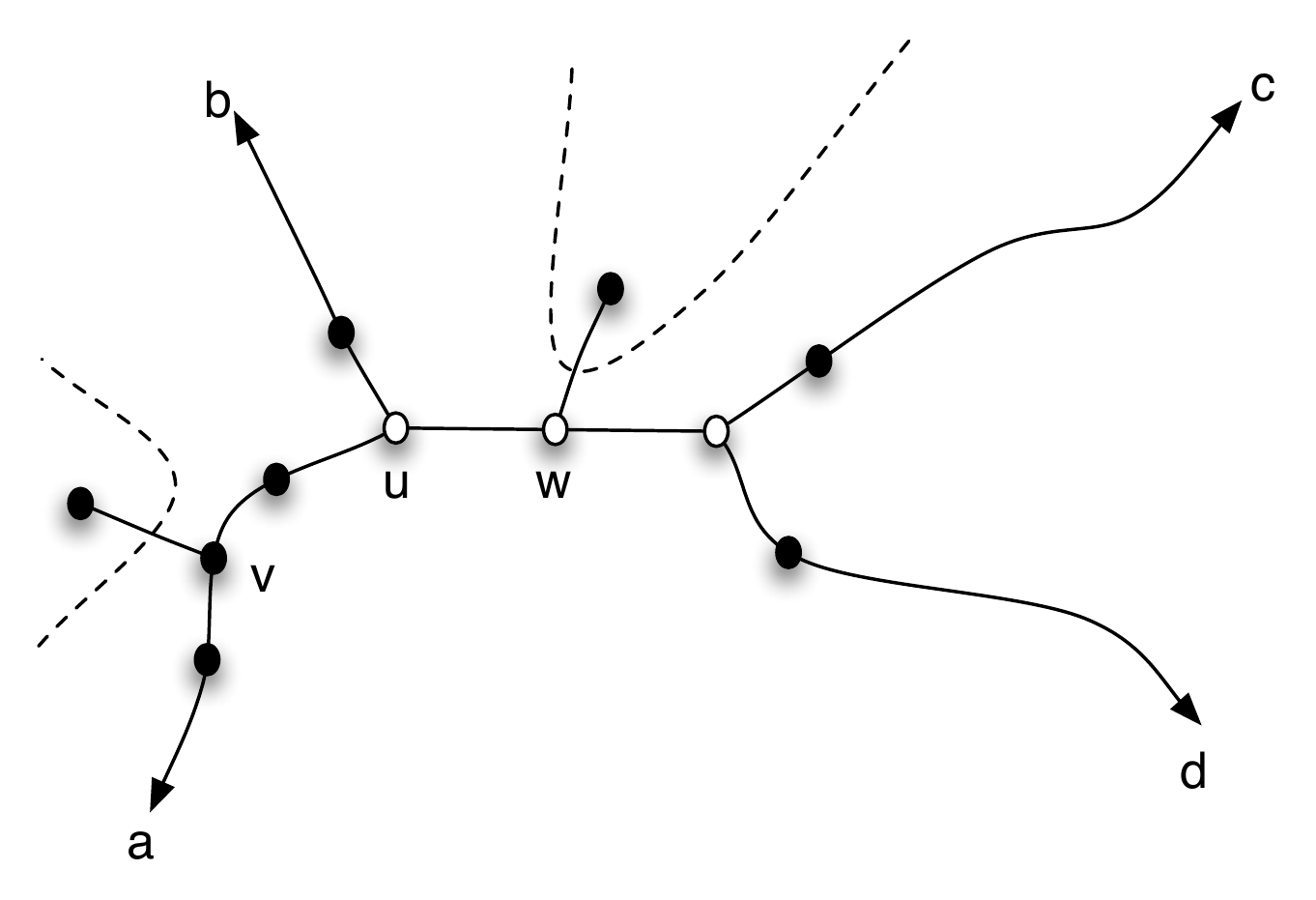}
\end{center}
\caption{Illustration for the three types of orbits of $\Aut({\mathcal M})_{(A)}$ which are disjoint from $A$, for $A = \{a,b,c,d\}$.} \label{fig:orbits} \end{figure}

\item[(v)] We shall freely use without explicit justification observations such as the following: given finite $A\subset M$ and $x\in M\setminus A$, there is a cone $U$ containing $x$ and disjoint from $A$; such a cone can be chosen properly inside any cone which contains $x$ and is disjoint from $A$, or, if $|A|\geq 3$, can be chosen to be at an $A$-vertex or at a vertex of $S(u,a)$ for some $A$-centre $u$ and $a\in A$.
\end{enumerate}
\item For any finite $A\subset M$, $\Aut(\calM)_{(A)}$ has no finite orbits on $M\setminus A$.
 This is immediate from  (2) $(iv)$. It follows that for any $A\subset M$, the model-theoretic algebraic closure of $A$ (the union, denoted $\acl(A)$, of the finite $A$-definable subsets of $M$) is exactly $A$.
 This property is inherited by reducts of $\calM$.

\item The structure $\calM$ is NIP, that is, its theory does not have the independence property
(see p. 69 of \cite{shelah} for  a definition). 
One way to see this is to observe that a model of $\Th(\calM)$  is interpretable in the field ${\mathbb Q}_2$, since it is well-known that the $p$-adic fields ${\mathbb Q}_p$ are  NIP.
Indeed, an elementarily equivalent structure $\calM'$ lives on the projective line ${\rm PG}_1({\mathbb Q}_2)$, with the $D$-relation defined by whether or not the cross-ratio lies in the maximal ideal (see Section~7 of \cite{m2}, or \cite[30.4]{an}). This observation is related to Proposition~\ref{valuation} below, and the remarks after it. Alternatively, as noted by the referee, observe that the set  $M$ in Example~\ref{mainex}(1) may be viewed as the set of leaves of a tree with internal vertices the set $M_0$ of all finite support functions from $(-\infty,n)$ (for some $n\in {\mathbb Z}$) to $\{0,1\}$, that such a `coloured' tree has NIP theory e.g. by \cite{parigot}, and that $(M,D)$ is interpretable in this tree.
\end{enumerate}
\end{remark}

A structure theory (or classification, in a loose sense) for primitive Jordan groups is given in \cite{am}, but in this paper  we just require the following consequence (Theorem~\ref{jordan} below), for which we first give some definitions.

A {\em separation relation} is the natural quaternary relation of separation on a circularly ordered set. An example  is the relation $S$ defined in the introduction. Following Section~3 of \cite{an}, a separation relation may be viewed as a relation satisfying the following universal axioms (though here, unlike the relation $S$ in the introduction, we do not require all arguments to be distinct).
\begin{enumerate}
\item[(S1)] $S(x,y;z,w) \to (S(y,x;z,w)\wedge S(z,w;x,y))$;
\item[(S2)] $(S(x,y;z,w) \wedge S(x,z;y,w))\leftrightarrow (y=z \vee x=w)$;
\item[(S3)] $S(x,y;z,w) \to (S(x,y;z,t)\vee S(x,y;w,t))$;
\item[(S4)] $S(x,y;z,w) \vee S(x,z;w,y) \vee S(x,w;y,z)$. 
\end{enumerate}

A {\em Steiner system} on $X$ is a Steiner $k$-system for some integer $k\geq 2$, that is, a family of subsets of $X$, called {\em blocks}, all of the same size (possibly infinite), such that any $k$ elements of $X$ lie on a unique block. It is {\em non-trivial} if any block has cardinality greater than $k$ and there is more than one  block. 

Next, following 
 \cite[Definition 2.1.10]{am}, we say that $H<\Sym(X)$ preserves a {\em limit of Steiner systems} on $X$ if for some $n>2$, $(H,X)$ is $n$-transitive but not $(n+1)$-transitive, and there is a totally ordered index set $(J,\leq)$ with no greatest element, and an increasing chain $(X_j:j\in J)$ of subsets of $X$ such that:
\begin{enumerate}
\item[(i)] $\bigcup(X_j:j\in J)=X$;
\item[(ii)] for each $j\in J$, $H_{\{X_j\}}$ is $(n-1)$-transitive on $X_j$ and preserves a non-trivial Steiner $(n-1)$-system on $X_j$;
\item[(iii)] if $i<j$ then $X_i$ is a subset of a block of the $H_{\{X_j\}}$-invariant Steiner $(n-1)$-system on $X_j$.
\item[(iv)] for all $h\in H$ there is $i_0\in J$, dependent on $h$, such that for every $i>i_0$ there is $j\in J$ such that $X_i^h=X_j$ and the image under $h$ of every block of the Steiner system on  $X_i$ is a  block of the Steiner system on $X_j$;
\item[(v)] for every $j\in J$, the set $X\setminus X_j$ is a Jordan set for $(H,X)$.
\end{enumerate}
Note that in \cite[Definition 2.1.10]{am} `$n$ is an integer greater than 3' should read `$n$ is an integer with $n\geq 3$'. 

Finally, we say a little more about general $D$-relations, as defined in the introduction. For a formal presentation, see Sections 22--26 of \cite{an}. Intuitively, there is a notion of `general betweenness relation' (see \cite[Sections 15-21]{an}) which can be defined on a (possibly dense) semilinear order with appropriate properties -- it is analogous to the well-known notion of $\Lambda$-tree (see \cite{chis}) but without any metric. A $D$-relation is the natural quaternary relation defined on a dense set of `ends', or `directions' (equivalence classes of 1-way paths) of the general betweenness relation. Following \cite[Section 24]{an}, if $D$ is a $D$-relation on a set $X$, then a {\em structural partition of $(X,D)$ with sectors $\{X_i:i\in I\}$} is a partition of $X$ into non-empty sets $X_i$ (with $|I|\geq 3$) such that
\begin{enumerate}
\item[(i)] for all $i\in I$ and $x,y\in X_i$ and  $z,w\not\in X_i$, we have $D(x,y;z,w)$,
\item[(ii)] if $x,y,z,w$ all lie in different sets $X_i$, then $D$ does not hold of any permutation of $x,y,z,w$.
\end{enumerate}
Extending our usage above, we use the word `cone' in place of `sector'. Though not explicitly stated, it is clear from \cite[Sections 23-26]{an}, that if $D$ and $D'$ are $D$-relations on $X$ with  the same set of cones, then $D=D'$. In fact, for distinct $x,y,z,w\in X$, $D(x,y;z,w)$ holds if and only if some cone contains $x,y$ and omits $z,w$: here, the direction $\Leftarrow$  follows from (i) above, and for $\Rightarrow$, by \cite[Theorem 24.2]{an} there is a structural partition $\{X_i:i\in I\}$ with $y,z,w$ in distinct cones, and by (4) at the start of \cite[Section 24]{an}, this and $D(x,y;z,w)$ implies that some $X_i$ contains $x,y$ and omits $z,w$.

A family ${\cal F}$ of subsets of $X$ is called {\em syzygetic} \cite[pp. 117-118, Section 34]{an} if for any $U,V\in {\cal F}$, if $U\setminus V$, $V\setminus U$ and $U\cap V$ are all non-empty then $U\cup V=X$. If $H\leq \Sym(X)$ then $U\subset X$ is {\em syzygetic} (for $H$) if $\{U^h:h\in H\}$ is syzygetic. By \cite[Corollary 25.2(2)]{an}, the family of cones of a $D$-relation is syzygetic, and hence any cone of a $D$-relation $(X,D)$ is syzygetic for $\Aut(X,D)$. Syzygetic sets provide a tool for recognising $D$-relations. For example, by 
\cite[Theorem 34.7]{an}, if $G$ is a 2-transitive permutation group on an infinite set $X$ and $U\subset X$ is syzygetic for $G$, then $G$ preserves on $X$ a $D$-relation or $C$-relation or general betweenness relation (as defined in \cite{an}). 

The following theorem follows from Theorem 1.0.2 of \cite{am}. 

\begin{theorem} \cite{am} \label{jordan} 
Let $G$ be a  Jordan permutation group on an infinite set $X$, and suppose that $G$ is 3-transitive but not highly transitive.
 Then $G$ preserves on $X$ a separation relation, a $D$-relation, a Steiner system, or a limit of Steiner systems.
\end{theorem}

\section{Proof of Theorem~\ref{reducts} -- group reducts}\label{sect:group-reducts}
 
 We here prove Theorem~\ref{reducts} (1), using the classification of 3-transitive Jordan groups
 (Theorem~\ref{jordan}).
 Put $G:=\Aut(\calM)$. We suppose  that $H$ is a closed proper subgroup of $\Sym(M)$ containing $G$, and must show that $H=G$. First observe that $H$ acts 3-transitively on $M$, and that each cone of $\calM$ is a proper Jordan set for $H$. Thus, by Theorem~\ref{jordan}, $H$ preserves a separation relation, a $D$-relation, a Steiner $k$-system for some $k\geq 3$, or a limit of Steiner systems. We consider these in turn, showing that the only possibility is that $H$ preserves the original $D$-relation of $\calM$, so $H=G$.

The group $G$ does not preserve any separation relation on $M$, and hence neither can $H$. 
For let $U$ be a cone of $\calM$ containing distinct elements $u,v$, and let $x,y,z$ be distinct elements of $M\setminus U$. By Remark~\ref{notes}~(1) there is $g\in G$ fixing  each of $x,y,z$ and with $(u,v)^g=(v,u)$. However, it is easily checked that no automorphism of a separation relation can fix three points and interchange two others.

The group $H$ cannot preserve a non-trivial Steiner $k$-system: for if $S$ is such a Steiner system with point set $M$, let $a_1,\ldots,a_k$ be distinct points, let $l$ be the block through
$a_1,\ldots,a_k$, let $b_k$ be a point not on $l$, and let $m$ be the block through $a_1,\ldots,a_{k-1},b_k$. Since blocks have more than $k$ points, there is $a_{k+1} \not\in \{a_1,\ldots,a_k\}$ on $l$, and $b_{k+1}\not\in \{a_1,\ldots,a_{k-1},b_k\}$ on $m$.
Now $l$ is the unique block containing $a_1,\ldots,a_{k-2},a_k,a_{k+1}$ and $m$ is the unique block containing $a_1,\ldots,a_{k-2},b_k,b_{k+1}$. These two blocks have common points
$a_1,\ldots,a_{k-1}$. Put $A:=\{a_1,\ldots,a_{k-2},a_k,a_{k+1},b_k,b_{k+1}\}$. Thus, $a_{k-1}$ lies in a finite orbit (in fact of size one) of $H_{(A)}$. This contradicts Remark~\ref{notes}~(3), as $H\supseteq G$ and $a_{k-1}\not\in A$.

Next, suppose that $H$ preserves a limit of Steiner systems on $M$, with the notation $M=\bigcup (X_j:j\in J)$ used in the last section. Though it is not explicit in \cite{am}, we may suppose (after replacing $J$ by a subset of the form $\{j\in J:j>j_0\}$ if necessary) that all the $X_i$ are infinite. Indeed, otherwise $M$ is a union of a sequence of finite sets whose complements are cofinite Jordan sets, in which case, by the main theorem in Section~4 of \cite{neumann}, $G$ is a group of automorphisms of a non-trivial Steiner system on $M$.

No set $X_j$ can contain a cone $U$. For  otherwise,
pick $a\in U$, and a finite set $A\subset X_j$, such that $a\not\in A$ and 
every automorphism of the Steiner system on $X_j$ which fixes $A$ pointwise fixes
$a$ (this is done as in the Steiner system argument in the last paragraph). Replacing $U$ by a subcone if necessary (cf. Remark~\ref{notes}(2)(v)) we may suppose that $U\cap A=\emptyset$. Now 
$G_{((M\setminus U) }$ has an element $h$ such that $a^h\neq a$  and
 $h\in H_{((M\setminus X_j)\cup A)}\leq H_{\{X_j\},(A)}$, which is a contradiction, since the Steiner system on $X_j$ is $H_{\{X_j\}}$-invariant. 
So for all $j\in J$, the set $M\setminus X_j$ is dense in  $M$ (in the sense given in  Section~\ref{sect:intro}). 

Since $M$ is countable and $(X_j:j\in J)$ is an increasing sequence of sets 
ordered by inclusion, $J$ has a countable cofinal subset $I=\{i_n:n\in \omega\}$. 
We claim that any infinite subset $A$ of $M$ meets infinitely many disjoint
 cones: indeed, the tree $T_0$ of $A$-vertices is infinite, so (e.g. by K\"onig's Lemma) contains an infinite path including infinitely many $A$-centres, and we may choose disjoint cones, one at each $A$-centre on this path, so that each meets $A$. 
In particular, $X_{i_0}$ meets infinitely many disjoint cones 
$\{U_i:i\in \omega\}$.  In particular, for each $n\in \omega$ there is
 $x_n\in U_n\cap X_{i_0}$. By the last paragraph there is
 $y_n\in U_n\setminus X_{i_n}$ for each $n\in \omega$. Since cones are Jordan sets and we may piece together actions on disjoint Jordan sets to obtain an automorphism, there is $g\in G$ such that $x_n^g=y_n$ for all $n$.
Thus, for each $n$, $X_{i_0}^g$ meets $M\setminus X_{i_n}$, so for each $j\in J$ with 
$j\geq i_0$ and each $k\in J$, $X_j^g \cap (M\setminus X_k)\neq \emptyset$. This contradicts clause 
$(iv)$ in the definition of a limit of Steiner systems.

Finally, suppose that $H$ preserves a $D$-relation $D'$ on $M$. We must show that $D'=D$. We may suppose $H=\Aut(M,D')$.

As noted at the end of Section 2, any $D$-relation is determined by its cones. Hence it suffices to show that $\calM=(M,D)$ and $(M,D')$ have the same cones. To see that every cone of $(M,D')$ is a cone of $(M,D)$, it suffices to show that if $V\subset M$ is syzygetic with respect to $G$ (see Section 2), and is infinite and coinfinite, then $V$ is a cone of $\calM$. Choose distinct $x,y\in V$, and $z\in M\setminus V$. Let 
$a:=\vertex(x,y,z)$ (with respect to $(M,D)$). Let $U_x,U_y,U_z$ be the cones of $\calM$ at $a$ containing $x,y,z$ respectively. Suppose first $U_x$ contains an element $w\not\in V$. By Remark~\ref{notes} (1), using that $M\setminus U_x$ is a cone of $\calM$, there is $g\in G_{(U_x)}$ such that $y^g=z$ and $z^g=y$. Then 
$x\in V \cap V^g$, $w\not\in (V \cup V^g)$, $y\in V\setminus V^g$ and $z\in V^g\setminus V$, contradicting that $V$ is syzygetic with respect to $G$. Thus, $U_x \subset V$, and similarly $U_y\subset V$. Now if $M\setminus V$ is not a cone of $(M,D)$, there is $w\neq z$ such that
$w\not\in V$, and such that if $U$ is the smallest cone of $(M,D)$ containing $w,z$, then $U$ also contains a point $t$ of $V$. We may suppose that $D(x,w;z,t)$ holds (otherwise $D(x,z;w,t)$ holds, and the same argument applies with $z$ and $w$ reversed).  Now let $U'$ be the smallest cone containing $z$ and $t$. There is $h\in G_{(U')}$   such that $x^h=w$ and $w^h=x$, and again, this contradicts that $V$ is syzygetic with respect to $G$. Thus, the conclusion is that $M\setminus V$ is a cone  of $(M,D)$, and hence also $V$ is a cone of $(M,D)$. 

Finally, as $\Aut(\calM)$ is transitive on the set of cones of $D$ (by Lemma~\ref{isom}(iv)) and preserves $D'$, and {\em some}
$(M,D)$-cone is an $(M,D')$--cone, it follows that {\em every} $(M,D)$-cone is an $(M ,D')$-cone. Thus $G=H$, and the proof of Theorem~\ref{reducts}~(1) is complete. $\Box$

\section{Proof of Theorem~\ref{reducts} -- definable reducts}
\label{sect:def-reducts}

Our proof below depends on the following analysis of definable subsets of $M$ in $\calM$.

\begin{lemma} \label{fraisse}
Let $A$ be a finite subset of $M$, and let $X$ be an infinite co-infinite $A$-definable subset of $M$ in the structure  $\calM$. Then there is an $A$-definable  subset $S$ of $M$ such that $S\triangle X\subseteq A$  and $S$ is a union of finitely many disjoint cones, each of which is  
\begin{enumerate}
\item[(i)] a cone, at an $A$-vertex, which is disjoint from $A$, or
\item[(ii)] for some $A$-centre $v$ and $a\in A$, a cone at a vertex of $S(v,a)$.
\end{enumerate}
\end{lemma}

Observe that any cone of type $(i)$ or $(ii)$ in the lemma is at a vertex which lies on a line of $T$ with both ends in $A$.
\begin{proof}
By Remark~\ref{notes}~(2(iv), 3), any infinite co-infinite  orbit of $G_{(A)}$ is 
\begin{enumerate}
\item[(a)] a cone at an $A$-vertex which is disjoint from $A$, or
\item[(b)] a cone, disjoint from $A$, at a vertex of $S(v,a)$ for some $A$-centre $v$ and $a\in A$. 
\end{enumerate}
Furthermore, an  Ehrenfeucht-Fra\"iss\'e game argument, which we omit, shows the following:

(*) for any $A$-centre $v$ and $a\in A$, there is a vertex $w\in S(v,a)$ such that if $U_w$ is the cone at $w$ containing $a$, then either $U_w\setminus \{a\}\subset X$ or $(U_w\setminus \{a\})\cap X=\emptyset$.

\noindent
We do not give full details of this game-theoretic argument, but it is standard. It suffices to observe the following: with $v,a$ as in (*), let $\bar{a}$ enumerate $A$, suppose that there is no $A$-vertex on $S(v,a)$ except $v$, that $m\in {\mathbb N}$, and that for $i=1,2$ there are $w_i\in S(v,a)$ with $d(v,w_i)\geq 2^m$, cones $U_{w_i}$ at $w_i$ disjoint from $A$, and $b_i\in U_{w_i}$; then 
$\bar{a}b_1\equiv_m \bar{a}b_2$, that is, $\bar{a}b_1$ and $\bar{a}b_2$ satisfy the same formulas of quantifier rank at most
 $m$ in the language $\{D\}$ of $\calM$. To start the argument, observe that two tuples $\bar{a}$, $\bar{b}$ of distinct elements of $M$ have the same quantifier-free type if and only if the trees that they `generate' (with vertices the $\bar{a}$-centres and $\bar{b}$-centres and edges the paths betwen these centres) are isomorphic as graphs. For an account of such arguments, see \cite[Section 3.3]{hodges}.

More informally, (*) asserts that if we consider cones of type (b) at vertices $w$ on $S(v,a)$, 
then either such cones all lie in $X$ provided $d(v,w)$ is sufficiently large, or they are all disjoint from $X$ for $d(v,w)$ sufficiently large.

We may assume that $X$ is disjoint from $A$.
Since $X$ is a union of $G_{(A)}$-orbits, it is a union of sets of type (a) or (b). Since there are finitely many $A$-vertices, there are finitely many cones of type (a), and these are of type $(i)$ in the lemma. There are also finitely many pairs of form $(v,a)$ with $v$ an $A$-centre and $a\in A$. For any such $(v,a)$, by (*), the union of the set of cones of type (b) can be written as the union of finitely many cones of type (b), together possibly with a set of form $W\setminus\{a\}$, where $W$ is a cone containing $a$ at a vertex of $S(v,a)$. Thus, the union of the sets of type (b) contained in $X$ can be written as the union of finitely many cones of type $(ii)$, possibly adjusted by the  removal of finitely many elements of $A$. This yields the lemma. 
\end{proof}

{\em Proof of Theorem~\ref{reducts}~(2). }
The proof proceeds in a series of claims.  We suppose that $\calM'$, a structure in a language $L'$,  is a non-trivial definable reduct of $\calM$. 
Our goal is to show that if $\calM'$ is not a trivial reduct (that is, a pure set), then the set of cones in $\calM$ is uniformly definable in $\calM'$. For then, by Claim~\ref{claim:cones-suffice} below, the relation $D$ is definable in $\calM'$, so the latter is an improper reduct.
 
\begin{claim}\label{claim:cones-suffice}
 Suppose that  the set of cones of $\calM$ is uniformly definable in $\calM'$; that is, for some integer $k>0$ there are $L'$-formulas $\phi(x,\bar{y})$ (with $\bar{y}=(y_1,\ldots,y_k)$) and $\psi(\bar{y})$, both over $\emptyset$,  such that the set
 $$\calF:=\{\phi(M,\bar{a}): \bar{a}\in M^k, \calM' \models \psi(\bar{a})\}$$
 is exactly the set of cones of $\calM$. 
 Then $D$ is definable in $\calM'$, so 
 $\calM'$ is an  improper definable reduct of $\calM$.
 \end{claim}
\begin{proof} 
 Define $D'(x,y;z,w)$ to hold on $M$ if one of
 \begin{enumerate}
\item[(a)] $x=y \wedge x\neq z \wedge x\neq w$, 
\item[(b)] $z=w\wedge z\neq x \wedge z\neq y$,
\item[(c)] $x,y,z,w$ are distinct and there are disjoint cones $U,U'$ with $x,y\in U$ and $z,w\in U'$.
 \end{enumerate}
 It is routine to check that $D'$ is exactly the relation $D$. 
 \end{proof}
 
 \begin{claim}\label{claim:cone-definable}
There is a cone of $\calM$ which is definable in $\calM'$.
\end{claim}
\begin{proof}
We first claim that in the structure $\calM'$ there is an infinite co-infinite definable subset of $M$. Indeed, by Remark~\ref{notes}~(3), for any $A\subset M$ the algebraic closure of $A$ in the structure $\calM'$ is exactly $A$. 
As $\calM'$ is a non-trivial definable reduct, there is finite $A\subset M$ enumerated by a tuple $\bar{a}$, an $L'$-formula $\phi(x,\bar{y})$, and $b,c\in M\setminus A$, such 
that $\calM'\models \phi(b,\bar{a})\leftrightarrow \neg \phi(c,\bar{a})$. 
Let $X:=\{x\in M: \calM'\models \phi(x,\bar{a})\}$. Then the sets $X\setminus A$ and $M\setminus (X\cup A)$ are both non-empty and definable, so as $\acl(A)=A$, these sets are both infinite. In particular, $X$ is an infinite co-infinite definable set in $\calM'$.

Thus, let $X$ be an infinite co-infinite subset of $M$ defined in $\calM'$ be the $L'$-formula $\rho(x,\bar{a})$.
By Lemma~\ref{fraisse}, replacing $X$ by a set differing finitely from it if necessary, we may write $X$ as a finite disjoint union $X=U_1\cup\cdots\cup U_r$ of cones, with each cone $U_i$ at the vertex $v_i$. Let $S:=\{v_1,\ldots,v_r\}$, and let $T_0$ be the finite subtree of $T$ consisting of the vertices of $S$ and the vertices and edges on paths between vertices in $S$. Define $s$ to be the sum of the distances (in the tree $T_0$) between distinct members of $S$.
We suppose that $\bar a$ and the definable set $X$ have been chosen to minimise $s$. We may also suppose that $s>0$; indeed, if $s=0$ then $X$ is the union of at most two cones at a vertex, so is a cone.

The finite tree $T_0$ has a leaf, $v_1$ say. Now since the $U_i$ are disjoint and their union is coinfinite in $M$, no vertex of $S$  lies in $U_1$. (Recall here the abuse of notation mentioned early in Section~\ref{sect:jordan}.) Furthermore, $v_j\neq v_1$ for $j>1$. For otherwise,
$U_1\cup U_j$ is a cone $W$ at a vertex $w$ (a neighbour of $v_1$); we may then replace $U_1\cup U_j$ by the cone $W$,
 contradicting the minimality of $s$.

It follows that there is a cone $U_1'$ at $v_1$ which is disjoint from $X$.
As in the argument at the end of Remark~\ref{notes}(1),  there is $g\in \Aut(\calM)$ fixing $M\setminus (U_1\cup U_1')$ pointwise, with $U_1^g=U_1'$ and $U_1'^g=U_1$.
 Let $X':=X\cup X^g= X\cup U_1'$, and put $\bar{a}':=\bar{a}^g$. Then $X'$ is $\bar{a}\bar{a}'$-definable, and it is easily checked that $X'$ is infinite and coinfinite. Furthermore, in the description of $X'$ as a disjoint union of cones we may replace $U_1$ by $W:=U_1\cup U_1'$ and thereby reduce $s$, contradicting our minimality assumption.
 \end{proof}

\begin{claim}\label{claim:clean-cone}
There is a cone of $\calM$ definable by an $L'$-formula $\phi(x,\bar{a})$ such that no element of $\bar{a}$ lies in $U$.
\end{claim} 
\begin{proof}
By Claim~\ref{claim:cone-definable}, there is a cone $U_1$ defined by $\phi(M,\bar{a})$ in $\calM'$. We may assume some element of $\bar{a}$ lies in $U_1$. Suppose $U_1$ is a cone at the vertex $u$ of $T$, and let $U_2$ and $U_3$ be the other two cones at $u$. Suppose first that one of these cones, say $U_3$, contains no element of $\bar{a}$. By for example the remarks after Lemma~\ref{isom}, there is $g\in G$ interchanging $U_1$ and $U_2$, and so fixing $U_3$ setwise: indeed, pick $b_1\in U_1$, $b_2\in U_2$ and $b_3\in U_3$, and choose $g$ so that
 $(b_1,b_2,b_3)^g=(b_2,b_1,b_3)$. Now $U_2$ is $\bar{a}^g$-definable, so $U_1 \cup U_2$ is 
 $\bar{a}\bar{a}^g$-definable, and hence $U_3=M\setminus(U_1 \cup U_2)$ is $\bar{a}\bar{a}^g$-definable -- and the parameter set $\bar{a}\bar{a}^g$ does not meet $U_3$. 
 
 Thus, we may suppose that the entries $a_1,a_2,a_3$ of $\bar{a}$ lie in $U_1,U_2,U_3$ respectively. Let $z$ be the vertex on the set $S(u,a_3)$  nearest to $u$ such that some cone $Z$ at $z$ contains no element of $\bar{a}$. We may suppose that $\phi(x,\bar{a})$ was chosen to minimise $d=d(u,z)$ (here we are allowing the formula $\phi$ to vary, but work with the above framework of formulas defining cones). Now let $u_1$ be the vertex adjacent to $u$ on the path from $u$ to $z$ (so possibly $u_1=z$). Let $h\in G_{(U_3)}$ with $(a_1,a_2)^h=(a_2,a_1)$ (such $h$ exists, for example by Remark~\ref{notes}(1)). Then $U_1\cup U_2$, a cone at $u_1$, is defined by the formula $\phi(x,\bar{a})\vee \phi(x,\bar{a}^h)$. The parameters $\bar{a}\bar{a}^h$ of this formula  have no entries lying in $U_3$ other than those of $\bar{a}$ which lie in $U_3$. Thus, we can replace the cone $U_1$ by $U_1 \cup U_2$, replace $u$ by $u_1$, and replace $\phi(x,\bar{a})$ by $\phi(x,\bar{a})\vee \phi(x,\bar{a}^h)$, and we have reduced $d$. This contradiction to minimality completes the proof of the claim. 
 \end{proof}
 
 Let $n:=l(\bar{a})$, where $\phi(x,\bar{a})$ is as in Claim 3. Since $G:=\Aut(\calM)$ is transitive on the set of cones, {\em every} cone $U$ of $\calM$ has the form $\phi(M,\bar{a}')$ for some $\bar{a}'\in (M\setminus U)^n$. However, for some $\bar{a}'$, the set $\phi(M,\bar{a}')$ might not be a cone, so we cannot immediately apply Claim~\ref{claim:cones-suffice} to define $D$. 
 Now let 
$\psi_1(\bar{y})$ be the $L'$-formula
 $$\exists z (\bigwedge_{i=1}^n z\neq y_i \wedge \neg \phi(z,\bar{y}))\wedge 	\exists z \, \big (\bigwedge_{i=1}^n z\neq y_i\wedge \phi(z,\bar{y}) \big ) \wedge \bigwedge_{i=1}^n \neg \phi(y_i,\bar{y}).$$
Let $\phi_1(x,\bar{y})$ be $\phi(x,\bar{y}) \wedge \psi_1(\bar{y})$. Now $\phi_1$ and $\psi_1$ are the first approximations of the formulas $\phi$ and $\psi$ mentioned in Claim~\ref{claim:cones-suffice}, and we shall repeatedly modify them until we obtain formulas as in Claim~\ref{claim:cones-suffice}.
 Observe by Remark~\ref{notes}~(3) that for any $\bar{y}$,  if $\phi_1(x,\bar{y})$ holds for some $x$ then $\phi_1(M,\bar{y})$ is infinite and coinfinite. So for any $\bar{a}'$, $\phi_1(M,\bar{a}')$ if non-empty is an infinite coinfinite set disjoint from $\bar{a}'$. Also, by transitivity of $\Aut(\calM)$ on the set of cones, {\em every} cone of $\calM$ has the form $\phi_1(M, \bar{a}')$ for some $\bar{a}'$ such that
 $\calM'\models \psi_1(\bar{a}')$.

Next let $\psi_2(\bar{y})$ be the following $L'$-formula, where $\bar{y}'=(y_1',\ldots,y_n')$:
$$ \forall \bar{y}'(\bigwedge_{i=1}^n \neg \phi_1(y_i',\bar{y}) \to ((\phi_1(M,\bar{y}')\supseteq \phi_1(M,\bar{y})) \vee (\phi_1(M,\bar{y}')\cap \phi_1(M,\bar{y})=\emptyset))).$$
Observe that $\psi_2(\bar{a})$ holds. For $\psi_1(\bar{a})$ holds, and if $\bar{a}'$ is disjoint from $\phi_1(M,\bar{a})$ then as $\phi_1(M,\bar{a})$ is a cone so a Jordan set,
$G_{\bar{a},\bar{a}'}$ is transitive  on it; hence if $b\in \phi_1(M,\bar{a})\cap \phi_1(M,\bar{a}')$ and $c\in \phi_1(M,\bar{a})$, there is $g\in G_{\bar{a},\bar{a}'}$ with $b^g=c$, so also $c\in \phi_1(M,\bar{a}')$.
Let $\phi_2(x,\bar{y})$ be $\phi_1(x,\bar{y}) \wedge \psi_2(\bar{y})$. Again, {\em every} cone of $\calM$ has the form $\phi_2(M, \bar{a}')$ for some $\bar{a}'$ such that
 $\calM'\models \psi_2(\bar{a}')$.

Let $\psi_3(\bar{y})$ be
$$\psi_2(\bar{y}) \wedge \exists \bar{z}\forall x \, \big(\phi_2(x,\bar{z})\leftrightarrow \neg \phi_2(x,\bar{y}) \big) \, .$$
Then since $\phi_2(M,\bar{a})$ is a cone, the complement of any cone is a cone, and $G$ is transitive on the set of cones, we have $\psi_3(\bar{a})$.
Let $\phi_3(x,\bar{y})$ be $\phi_2(x,\bar{y})\wedge \psi_3(\bar{y})$ (an $L'$-formula). Then we have
$$\psi_3(\bar{y})\to \exists\bar{z}(\forall x(\phi_3(x,\bar{z})\leftrightarrow \neg\phi_3(x,\bar{y}))).$$
Also,  by transitivity on the set of cones, every cone has the form $\phi_3(M,\bar{a}')$ for some $\bar{a}'$ such that $\calM'\models \psi_3(\bar{a}')$.

 \medskip
 
 We aim next to reduce to the case when $n=|\bar{y}|=4$. For this, it suffices to show that {\em some} cone $U$ is {\em 4-definable}, that is, definable by an $L'$-formula with 4 parameters, none lying in $U$; for then we may take this formula to be $\phi(x,\bar{y})$ above, and modify it to obtain
 $\phi_3(x,\bar{y})$ and $\psi_3(\bar{y})$ as above, but with $l(\bar{y})=4$.
 
 \begin{claim}\label{claim:well-definable} 
 Let $b_1,b_2,b_3$ be distinct elements of $M$, and for each $i$ let
$U_i$ be the cone at $u:=\vertex(b_1,b_2,b_3)$ containing $b_i$. 
Let $U_1$ be definable over $b_1,b_2,b_3$, with $U_1=\chi(M,b_1,b_2,b_3)$ where $\chi$ is an $L'$-formula. Then some cone is 4-definable.

\end{claim}
\begin{proof}
Choose $b_4$ so that if $u':=\vertex(b_1,b_3,b_4)$ then $u'\in S(u,b_3)$ and $d(u,u')=2$. Let $U_3',U_4'$ be respectively the cones 
at $u'$ containing $b_3,b_4$. Then, by 3-transitivity of $\Aut(\calM)$ each of the cones $U_1,U_2,U_3',U_4'$ is definable over $b_1,b_2,b_3,b_4$, respectively by the formulas $\chi(x,b_1,b_2,b_3)$, $\chi(x,b_2,b_1,b_3)$, $\chi(x, b_3,b_1,b_4)$, and $\chi(x, b_4, b_1,b_3)$. Hence $W:=M\setminus (U_1 \cup U_2\cup U_3'\cup U_4')$ is a cone (based at the common neighbour of $u$ and $u'$) which is $L'$-definable over $b_1,b_2,b_3,b_4$ but does not contain any of these parameters.
\end{proof}

Let $\phi_3$ and $\bar{a} \in M^n$ be as above, 
so that $\phi_3(M,\bar{a})$ is a cone. 
Now define the relation $E(x,y;z,w)$ to hold if and only if $x,y,z,w$ are distinct and there is $\bar{e} \in M^n$ such that $\phi_3(M,\bar{e})$ contains $x,y$ but not $z,w$. Observe that $E$ is $L'$-definable.
\begin{claim}\label{claim:distance}
 Suppose that $E(b_1,c;b_2,b_3)\wedge D(b_1,b_3;c,b_2)$ holds. 
Then $u:=\vertex(b_1,b_3,c)$ and $v:=\vertex(b_1,b_2,c)$ are at  distance at most $2n+1$.
\end{claim}
\begin{proof}
Suppose  for a contradiction that  $d(u,v)=m>2n+1$. There is $\bar{e}\in M^n$ such that
$\phi_3(M,\bar{e})$ contains $b_1,c$ and omits $b_2,b_3$. There are at most $n$ vertices strictly between $u$ and $v$ of form $\vertex(b_1,b_2,e_i)$ for some $e_i$ in $\bar{e}$. In particular, there is a vertex $w$ strictly between $u$ and $v$ and not of this form (in fact, there are at least $n+1$ such $w$). Let $W$ be the cone at $w$ which does not contain $b_1,b_2,b_3$ or any element of $\bar{e}$, and for $i=1,2$ let $W_i$ be the cone at $w$ containing $b_i$.  We claim that $W\subseteq \phi_3(M,\bar{e})$.
Indeed, suppose not. Then as $\Aut(\calM)_{\bar{e}}$ is transitive on $W$, we have $W\cap \phi_3(M,\bar{e})=\emptyset$. We may (reordering if necessary) write $\bar{e}=\bar{e}_1\bar{e}_2$ so that for $e_i\in \bar{e}_1$, the vertex $w$ lies on the $e_ib_2$ line, and for $e_i \in \bar{e}_2$, the element $w$ 
lies on the $e_ib_1$ line, and we may suppose $\bar{e}_2$ is non-empty (for clearly $\bar{e}\neq \emptyset$). Using Remark~\ref{notes}(1), pick $g\in G$
 fixing pointwise the cone $W_1$, and  interchanging $W$ and $W_2$. Put $\bar{e}':=\bar{e}^g$. Then $\bar{e}'=\bar{e}_1\bar{e}_2^g$ is disjoint 
from $\phi_3(M,\bar{e})$, as  $\bar{e}_1$ is disjoint from $\phi_3(M,\bar{e})$ and $\bar{e}_2^g$ lies in $W$ which is also disjoint from $\phi_3(M,\bar{e})$. However
 $\phi_3(M,\bar{e}')$ neither contains nor is disjoint from $\phi_3(M,\bar{e})$, as $b_1\in \phi_3(M,\bar{e}')\cap \phi_3(M,\bar{e})$, and $c\in \phi_3(M,\bar{e})\setminus \phi_3(M,\bar{e}')$. This contradicts  the assumption that  $\psi_2(\bar{e})$ holds, so yields that $W\subseteq \phi_3(M,\bar{e})$.

Since $\psi_3(\bar{e})$ holds, there is $\bar{e}'$ such that $\neg\phi_3(M,\bar{e})=\phi_3(M,\bar{e}')$. The set $\phi_3(M,\bar{e}')$ contains $b_2,b_3$ and omits $c,b_1$. Since there are at least $2n+1$ vertices strictly between $u$ and $v$, and $|\bar{e}|+|\bar{e}'|=2n$, there is a vertex $z$ strictly between $u$ and $v$ such that the cone $Z$ at $z$ omitting all the $b_i$  does not meet $\bar{e}$ or $\bar{e}'$. The argument of the last paragraph now shows that $Z\subseteq \phi_3(M,\bar{e})\cap \phi_3(M,\bar{e}')$, which is impossible. 
\end{proof}

\begin{claim} \label{claim:well-definable-two}
Some cone is 4-definable.
\end{claim}
\begin{proof} Let $b_1,b_2,b_3$ be distinct, let $u:=\vertex(b_1,b_2,b_3)$, and let $U_i$ be the cone at $u$ containing $b_i$, for each $i=1,2,3$. Put $X:=E(b_1,M;b_2,b_3) \cup\{b_1\}$. 
For any $c\in U_1\setminus\{b_1\}$, there is a cone containing $b_1,c$ and omitting $b_2,b_3$. Thus, by the definition of $E$ and the fact that every cone has form $\phi_3(M,\bar{a}')$ for some $\bar{a}'$, we have 
$X \supseteq U_1$. By Claim~\ref{claim:well-definable}, we may suppose $X\neq U_1$. Thus, we may suppose that $X\cap U_3\neq \emptyset$. By Claim~\ref{claim:distance} together with Lemma~\ref{fraisse},  $X\cap U_3$ is a union of some cones (not containing  $b_3$) at vertices of $S(u,b_3)$ which are at  distance at most $2n+1$ from $u$. 
Pick a vertex $w$ on $S(u,b_3)$ with $d(u,w)$ minimal such that the cone $W$ at $w$ which omits all the $b_i$ is disjoint from $X$.
 Using Remark~\ref{notes}(1), let $g\in \Aut(\calM)_{(U_3)}$ with $(b_1,b_2)^g=(b_2,b_1)$ and let $Y:=X \cup X^g$. Then $U_1 \cup U_2 \subset Y$, and $Y \cap W=\emptyset$; in fact, the cone at $w$ containing $b_1$ lies in $Y$, since by the choice of $W$, for any vertex $w'$ strictly between $u$ and $w$, the cone $W'$ at $w'$ not containing any $b_i$ has a non-empty intersection with $X$, so by Lemma~\ref{fraisse}, $W'$ even lies in $X$.
 
 Choose $b_4\in M \setminus \{b_1,b_2,b_3\}$ such that if $u':=\vertex(b_1,b_3,b_4)$ then $u'\in S(w,b_3)$ and $d(w,u')=d(w,u)$. Let $h\in G$ with $(b_1,b_2,b_3,b_4)^h=(b_4,b_3,b_2,b_1)$. Observe that $h$ fixes $W$ setwise (and in fact, can be chosen to fix $W$ pointwise).
Then as $Y$ is $b_1b_2b_3$-definable, the set $Y \cup Y^h$ is $b_1b_2b_3b_4$-definable. The claim follows, as $W:=M\setminus (Y \cup Y^h)$. 
\end{proof}

By Claim~\ref{claim:well-definable-two}, we now suppose that in the formula $\phi_3(x,\bar{y})$, the tuple $\bar{y}$ has length 4.
We shall show that either $\phi_3$ uniformly defines the family of cones, or it fails to do so in a very special way, and can be modified to give a formula which uniformly defines the set of cones.

\begin{claim}\label{claim:special}
Suppose that there is $\bar{b}\in M^4$ such that $\phi_3(M,\bar{b})$ is non-empty and is not a cone. Then there are adjacent vertices $w_1,w_2$ and disjoint cones $W_1,W_2$ at $w_1,w_2$ such that $\phi_3(M,\bar{b})=W_1 \cup W_2$. Furthermore (re-ordering $\bar{b}$ if necessary), we may suppose $D(b_1,b_2;b_3,b_4)$ holds, and $w_1,w_2$ lie  on the path between $\vertex(b_1,b_2,b_3)$ and $\vertex(b_1,b_3,b_4)$ which are at distance 3.
\end{claim}
\begin{proof} 
We may suppose that $D(b_1,b_2;b_3,b_4)$ holds. Let $u:=\vertex(b_1,b_2,b_3)$ and $v:=\vertex(b_1,b_3,b_4)$, and let $L$ be the path in $T$ between $u$ and $v$. Also let $V_1,V_2$ be the cones at $u$ containing $b_1,b_2$ respectively, and $V_3,V_4$ be the cones at $v$ containing $b_3,b_4$ respectively.

For each $i=1,\ldots,4$, there is no cone $U$ containing $b_i$  such that
$U\setminus \{b_i\}\subseteq\phi_3(M,\bar{b})$. For otherwise, as $\calM'\models \psi_3(\bar{b})$, there would be $\bar{c}$ such that
$M\setminus\phi_3(M,\bar{b})=\phi_3(M,\bar{c})$. As $\calM'\models \neg \phi_3(b_i,\bar{b})$ (since $\psi_1(\bar{b})$ holds), we would have $\phi_3(b_i,\bar{c})$, so $b_i\not\in \bar{c}$ (again as $\psi_1(\bar{c})$). Let $V$ be any subcone  of $U$ containing $b_i$ and disjoint from $\bar{c}$. Then $\phi_3(M,\bar{c})\cap V=\{b_i\}$. As $V$ is a Jordan set, this
 is impossible: indeed, if $d\in V\setminus\{b_i\}$ pick $g\in \Aut(\calM)_{(M\setminus V)}$ with $b_i^g=d$; then $\bar{c}^g=\bar{c}$, so $\phi_3(d,\bar{c})$ holds, a contradiction. 
 
It follows from the last paragraph and Lemma~\ref{fraisse} that there is $t$ such that $\phi_3(M,\bar{b})$ is a union of disjoint cones at vertices at distance at most $t$ from $u$ or $v$. (We do not here claim that $t$ is independent of $\bar{b}$; the existence of $t$ follows from the finiteness of the number of cones  in Lemma~\ref{fraisse}, and the fact that if $A=\{b_1,b_2,b_3,b_4\}$ then the only $A$-centres are $u$ and $v$.)
The argument now falls into two subcases.

{\em Subcase $(i)$.} Suppose that $\phi_3(M,\bar{b})$ contains a subcone $W$ of some $V_i$, say of $V_3$. 
Let $w:=\vertex(W)$. We may suppose that $w\in S(v,b_3)$, and that $W$ is chosen to minimise $d(v,w)$. 
 
Suppose first that  $\phi_3(M,\bar{b})\subset V_3$. In this case, as $\phi_3(M,\bar{b})$ is 
not a cone, there is some $w'\neq w$ on $S(v,b_3)$, and a cone $W'\subset \phi_3(M,\bar{b})$, with $\vertex(W')=w'$. Note that by minimality of $d(v,w)$, $w'\in S(w,b_3)$. We may choose 
$g\in G_{b_1b_3}$ (a translation along the $b_1b_3$ line) with $W'^g=W$. Then $b_2^g,b_4^g\not\in \phi_3(M,\bar{b})$, and $\phi_3(M,b_1b_2^gb_3b_4^g)$ meets $\phi_3(M,\bar{b})$ in a proper non-empty subset on $\phi_3(M,\bar{b})$. This contradicts that $\psi_2(\bar{b})$ holds. 

Thus, we may suppose $\phi_3(M,\bar{b})\not\subset V_3$.
Choose $z$ on $S(v,b_3)$ nearest to $v$ so that if $U$ is the cone at $z$ containing 
$b_1$, then $\phi_3(M,\bar{b})\subset U$. Let $Z$ be the cone at $z$ omitting $b_1$ and $b_3$ (so $Z\cap \phi_3(M,\bar{b})=\emptyset$), and pick $b_5\in Z$. 
Let $z'$ be the neighbour of $z$ on the $vz$-path, and let $Z'$ be the cone at $z'$ not containing  $b_1$ or $b_3$, so 
$Z'\subset \phi_3(M,\bar{b})$ (by minimality of $d(v,z)$). Choose a vertex $z''$ between $v$ and $z'$, as near to $z'$ as possible such that the cone $Z''$ at $z''$ not containing $b_1$ or $b_3$ does not lie in $\phi_3(M,\bar{b})$, and choose $b_6 \in Z''$ with $b_6\not\in \phi_3(M,\bar{b})$. (Possibly, $z''=v$ and $Z''=V_4$, in which case we choose $b_6=b_4$.) Observe that $D(b_1,b_6;b_3,b_5)$ holds. Let $g\in G$ with $(b_1,b_6,b_3,b_5)^g=(b_3,b_5,b_1,b_6)$, so $g$ interchanges  $z$ and $z''$.
Then $\bar{b}^g$ is disjoint from $\phi_3(M,\bar{b})$, and
$Z'\subset \phi_3(M,\bar{b})\cap \phi_3(M,\bar{b}^g)$. (Indeed,  if $z'''$ is the vertex between $z''$ and $z$ adjacent to $z''$, and $Z'''$ is the cone at $z'''$ not containing $b_1$ or $b_3$, then $Z'''\subset \phi_3(M,\bar{b})$ (by choice of $z''$) and $Z'=(Z''')^g\subset \phi_3(M,\bar{b})\cap \phi_3(M,\bar{b}^g)$.) However, $\phi_3(M,\bar{b}^g)$ lies in the cone at $z''$ containing $b_3$, so is contained in $V_3$ so disjoint from the non-empty set $\phi_3(M,\bar{b})\setminus V_3$. This contradicts that $\psi_2(\bar{b})$ holds.

{\em Subcase $(ii)$.}
Suppose that $\phi_3(M,\bar{b})$ is a union of cones at vertices of $L$. We consider first the case when there is a vertex $w\not\in\{u,v\}$ but lying on $L$, such that $\phi_3(M,\bar{b})$ contains no cone at $w$. In this case, let $W$ be the cone at $w$ not containing any $b_i$, so $\phi_3(M,\bar{b})\cap W =\emptyset$. We may suppose that $\phi_3(M,\bar{b})$ contains a cone at a vertex strictly between $w$ and $v$. It follows that $\phi_3(M,\bar{b})$ does not contain any cone $W'$ at a vertex between $w$ and $u$. For otherwise, pick  $g\in G_{b_3,b_4}$ with $b_1^g,b_2^g\in W$ and let $\bar{b}':=(b_1^g,b_2^g,b_3,b_4)$. Then $b_1^g,b_2^g\not\in \phi_3(M,\bar{b})$ and $\phi_3(M,\bar{b}')\cap\phi_3(M,\bar{b})$ is a proper non-empty subset of $\phi_3(M,\bar{b})$ (proper as it omits $W'$), contradicting that $\psi_2(\bar{b})$ holds. This argument, and appropriate choice of $w$ (to minimise $d(w,v)$),  reduces us to the case when there is $w'$ either equal to $v$ or strictly between $w$ and $v$
(on the line between them), such that
 $\phi_3(M,\bar{b})$ is exactly the union of the cones (by assumption, more than one of them), not containing any $b_i$,  at vertices strictly between $w$ and $w'$. In this case, let $u''$ be the neighbour of $u$ on $L$ and $v''$ be the neighbour of $v$ on $S(v,b_3)$. Choose $b_2''$ in the cone at $u''$ which does not meet $\bar{b}$, choose $b_4''$ in the cone at $v''$ which does not meet $\bar{b}$, and let $\bar{b}'':=(b_1,b_2'',b_4'',b_3)$. There is  $h\in G$ with $\bar{b}^h=\bar{b}''$, and for such $h$ we have $u^h=u''$ and $v^h=v''$ (in the induced action of $h$ on $T$). Again, we find that $\phi_3(M,\bar{b})\cap \phi_3(M,\bar{b}'')$ is a proper non-empty subset of $\phi_3(M,\bar{b})$. As $\bar{b}''$ is disjoint from $\phi_3(M,\bar{b})$ (by its choice), this contradicts that $\psi_2(\bar{b})$ holds.

Thus,  $\phi_3(M,\bar{b})$ is the union of all cones at vertices of $L$ which do not contain any $b_i$. In particular, we have shown that any non-empty set of form $\phi_3(M,\bar{c})$ which is not a cone has this form.

To prove the claim, we must now show $d(u,v)=3$. 
As $\phi_3(M,\bar{b})$ is not a cone, $d(u,v)\geq 3$ and the complement $M\setminus \phi_3(M,\bar{b})$ is a union of two cones at distinct vertices. As $\psi_3(\bar{b})$ hold, there is $\bar{b}'\in M^4$ such that 
$M\setminus \phi_3(M,\bar{b})=\phi_3(M,\bar{b}')$. However, if $d(u,v)>3$, then the set $\phi_3(M,\bar{b}')$ cannot have the forms described in the last paragraph. 
\end{proof}

We may now suppose  that there is $\bar{b}\in M^4$ such that $\phi_3(M,\bar{b})$ is a union of two disjoint cones $U_1,U_2$ at adjacent vertices $u_1,u_2$ respectively, as described in Claim~\ref{claim:special}; indeed, if there is no such $\bar{b}$, then by Claim~\ref{claim:special} the formula $\phi_3$ (with $\psi_3$) uniformly defines the set of cones in $\calM'$, and then Claim~\ref{claim:cones-suffice} completes the proof. 
We claim that $\Aut(\calM)$ is transitive on the collection ${\cal S}$ of all sets of the form $U\cup V$ where $U,V$ are disjoint cones at adjacent nodes with $U\cup V\neq M$; indeed, for two such pairs $(U,V)$ and $(U',V')$, there is $g\in \Aut(\calM)$ with $U^g=U'$, and $h\in \Aut(\calM)_{(U')}$ with $(V^g)^h=V'$, and $(U,V)^{gh}=(U',V')$.  Thus, every such set (in ${\cal S}$) has the form $\phi_3(M,\bar{b}')$ for some $\bar{b}'$.

Define $\chi(\bar{y})$ to be the $L'$-formula (with $l(\bar{y})=l(\bar{z})=4$)
\begin{align*}
\psi_3(\bar{y}) \wedge \exists \bar{z} \big ( & \exists u\phi_3(u,\bar{z}) \wedge \forall x(\phi_3(x,\bar{y}) \to \neg \phi_3(x,\bar{z})) \\
\wedge & \exists^{\leq 1}  i \; \phi_3(y_i,\bar{z}) \wedge \exists \bar{w}\forall x(\phi_3(x,\bar{w})\leftrightarrow (\phi_3(x,\bar{y}) \vee \phi_3(x,\bar{z}))) \big ) \, .
\end{align*}
The second conjunct says  there is $\bar{z}$ so that the non-empty set $\phi_3(M,\bar{z})$ is disjoint from $\phi_3(M,\bar{y})$ and contains at most one of $y_1,y_2,y_3,y_4$ and $\phi_3(M,\bar{y})\cup \phi_3(M,\bar{z})$ is in the family of sets defined by $\phi_3$ (so  is a cone or the union of two disjoint cones at adjacent vertices).

It can be checked that if $\phi_3(M,\bar{c})$ is a cone then $\chi(\bar{c})$ holds. Indeed, if $\phi_3(M,\bar{c})$ is the cone $U$ at $u$, then as $l(\bar{c})=4$ there is a vertex $v$ adjacent to $u$ in $T$ and not lying in $U$, and a cone $V$ at $v$ disjoint from $U$ and containing at most one element of $\bar{c}$. Let $V=\phi_3(M,\bar{c}')$. Then
$\bar{c}'$ is a witness for $\bar{z}$ in the second conjunct of $\chi$. Also $U\cup V$ is the union (not equal to $M$) of two disjoint cones at adjacent vertices so has the form $\phi_3(M,\bar{c}'')$ for some $\bar{c}''$ as required; here $\bar{c}''$ is a witness for $w$ in $\chi(\bar{c})$.

However, it can be checked that for $\bar{b}$ as above, with $U_1,U_2$ cones at adjacent vertices $u_1,u_2$ as above,
 $\chi(\bar{b})$ does not hold. Indeed, consider any non-empty set of the form $\phi_3(M,\bar{b}')$ which is disjoint from $\phi_3(M,\bar{b})=U_1\cup U_2$, and contains at most one element of $\bar{b}$. Then $\phi_3(M,\bar{b}')$ is a cone $V$ at a vertex $v$ or the union of two disjoint cones $V_1,V_2$ at adjacent vertices $v_1,v_2$. By the second assertion of Claim~\ref{claim:special}, the vertex $v$ (or $v_1,v_2$) cannot lie in $\{u_1,u_2\}$. It can now be  seen that $\phi_3(M,\bar{b}) \cup \phi_3(M,\bar{b}')$ is not a cone or the union of two disjoint cones at adjacent vertices.

Thus, the $L'$-formula $\phi_3(x,\bar{y})\wedge \chi(\bar{y})$ uniformly defines the family of all cones; that is, in Claim~\ref{claim:cones-suffice} we may take $\phi(x,\bar{y})$ to be $\phi_3(x,\bar{y})\wedge \chi(\bar{y})$, and $\psi(\bar{y})$ to be $\psi_3(\bar{y})\wedge \chi(\bar{y})$. By Claim~\ref{claim:cones-suffice}, this  completes the proof of the theorem.
$\Box$

\begin{remark} \label{groupsenough}\rm
1. We expect that the valency 3 assumption on $(T,R)$ is not needed, and that the $D$-relation arising from any regular combinatorial tree has no proper non-trivial definable reducts, but have not checked this. 
The assumption in Theorem~\ref{reducts} that the tree $(T,R)$ has degree 3 is used in various places in the proof, since we often use the facts that the complement of a cone is a cone, and the union of two cones at a vertex is a cone. It can be checked that if $(T,R)$ is any regular tree of degree at most $\aleph_0$ and at least three, $M$ is a dense set of ends of $T$, and $D$ is the induced $D$-relation on $M$, then $(M,D)$ has no proper non-trivial {\em group-reducts}. For the proof, a small adjustment is needed in the proof of Theorem~\ref{reducts} (1) in the case when a closed supergroup $H$ of $\Aut(M)$ preserves a $D$-relation $D'$; we omit the details.

2. It may be possible to obtain part (2) of Theorem~\ref{reducts} from a version of the proof of (1), done in a saturated model of ${\rm Th}(\calM)$, arguing as in Proposition~\ref{posetemb} below. (Such a countable saturated model exists and has analogous symmetry properties, with the two-way lines in the underlying betweenness relation carrying the linear betweenness relation induced from a countable saturated model of $({\mathbb Z},<)$.) Indeed, if $\calM$ had a proper non-trivial definable reduct, this would hold in any elementary extension. The automorphism group of a saturated elementary extension $\calM'$ of $\calM$ would also be a 3-transitive Jordan group, and would have (as in the proof of (1)) to preserve a $D$-relation. It would appear that the only `other' possible $D$-relation, apart from the natural one,  would be the one obtained by identifying elements of the underlying betweenness relation (definable in $\calM'$  as in \cite[Theorem 25.3]{an}) which are finitely far apart; this $D$-relation on $\calM'$ is not definable. Some further work is needed, and we view our proof of (2) above, not dependent on the very intricate structure theory for primitive Jordan groups, as of independent interest.

\end{remark}

\section{Other possible examples}
\label{sect:examples}
We discuss here some further approaches to constructing infinite non-$\omega$-categorical structures with no proper non-trivial reducts, focussing particularly on strongly minimal structures (defined in the introduction). The main result here is Theorem~\ref{dist-trans}.

In the next proposition,  the collection of definable reducts of a structure $\calM$ is partially ordered by putting $\calM_1\leq \calM_2$ if $\calM_2$ is $\emptyset$-definable in $\calM_1$; the group-reducts of $\calM$ are partially ordered by inclusion of groups.  We remark that if $\calM$ is strongly minimal over a countable language, then the first-order theory $T$ of $\calM$ is categorical in all uncountable cardinalities (Corollary 5.7.9. in~\cite{tz}) and, in particular, for every infinite cardinal $\kappa$, $T$ has a saturated model $N$ of cardinality $\kappa$.

\begin{proposition} \label{posetemb}
Let $\calM$ be a saturated structure. Then
\begin{enumerate}
\item[(i)] the partial order of definable reducts of $\calM$ embeds into the partial order of group-reducts of $\calM$, and
\item[(ii)] if $\calM$ has no proper non-trivial group-reducts, then $\calM$ has no proper non-trivial definable reducts.
\end{enumerate}
\end{proposition}

\begin{proof}
$(i)$ First observe that any definable reduct of a saturated structure is saturated. 
We show that if $\calM_1$ and $\calM_2$ are distinct definable reducts of $\calM$, then $\Aut(\calM_1)\neq \Aut(\calM_2)$, and  hence that if $\calM_2$ is a proper definable reduct of $\calM_1$ (both reducts of $\calM$) then $\Aut(\calM_1) < \Aut(\calM_2)$. 
To see these, let $k$ be least such that there is some $k$-ary relation $R$ which is $\emptyset$-definable in $\calM_1$ but not in $\calM_2$. 
Then by compactness and saturation, there are $\bar{a},\bar{b}\in M^k$ such that
 $\calM_1\models R\bar{a}\wedge \neg R\bar{b}$ and
 $\tp_{\calM_2}(\bar{a})=\tp_{\calM_2} (\bar{b})$.
 It follows by saturation that $\bar{a}$ and $\bar{b}$ lie in the same orbit of $\Aut(\calM_2)$ but in distinct orbits of $\Aut(\calM_1)$.

$(ii)$ This is immediate from $(i)$. 
\end{proof}

\begin{remark}\rm  Part (2) of Theorem~\ref{reducts} cannot be deduced from part (1) and Proposition~\ref{posetemb}, since the structure $\calM=(M,D)$ is not saturated, or even recursively saturated; however, see Remark~\ref{groupsenough}(2). The structure $(M,D)$ is homogeneous in the sense that any elementary map between finite subsets of $M$ extends to an automorphism, but this property might not {\em a priori} hold in  definable reducts.
\end{remark}

In this paper, we have considered definable reducts of a structure $\calM$ up to  interdefinability over $\emptyset$, and group-reducts, with two group-reducts identified if they have the same automorphism group. However, it is also possible to mix the two notions, and consider group-reducts of $\calM$ (structures $\calM'$ with the same domain as $\calM$ such that $\Aut(\calM') \supseteq \Aut(\calM)$) with two group reducts $\calM_1$ and $\calM_2$ identified if each is $\emptyset$-definable in the other; that is, group-reducts considered up to interdefinability. Under this notion,  the answer to the question of Junker and Ziegler mentioned in the introduction is positive.

\begin{proposition} \label{mix}
Let $\calM$ be a countable structure over a countable language and suppose that $\calM$  is not $\omega$-categorical. Then $\calM$ has $2^{\aleph_0}$ distinct group-reducts up to interdefinability.
\end{proposition}

\begin{proof}
Let $k$ be least such that $\Aut(\calM)$ has infinitely many orbits on $M^k$, and let these orbits be $\{P_i:i\in \omega\}$. For 
each $S\subset \omega$ let $\calM_S$ be the structure with domain $M$ and with a single $k$-ary relation 
interpreted by $\bigcup_{i\in S} P_i$. Then 
$\Aut(\calM)\leq \Aut(\calM_S)$, so $\calM_S$ is a group-reduct of $\calM$. As 
there are $2^{\aleph_0}$ such group-reducts and only countably many formulas in each of the languages, the result follows. 
\end{proof}

Next, we make an elementary observation on maximal-closed subgroups of $\Sym({\mathbb N})$. 

\begin{lemma}
Let $G$ be a maximal-closed proper subgroup of $\Sym({\mathbb N})$ which acts imprimitively on ${\mathbb N}$. Then $G$ acts oligomorphically on ${\mathbb N}$.
\end{lemma}

\begin{proof} If $G$ is intransitive on  ${\mathbb N}$ then as $G$ is maximal-closed, it is the setwise  stabiliser of a proper subset $A$ of ${\mathbb N}$, and $A$ must be finite or cofinite. 
Likewise, if $G$ is transitive but imprimitive on  ${\mathbb N}$ then $G$ is the stabiliser of a partition of ${\mathbb N}$ into parts of the same size, so is a wreath product of symmetric groups. In either case, it is routine to check (directly, or via the Ryll-Nardzewski Theorem) that $G$ acts oligomorphically on ${\mathbb N}$.
\end{proof}

We turn now to strongly minimal examples.
The goal is to find examples  of strongly minimal sets which are not $\omega$-categorical but
have no proper non-trivial reducts (in either of our two senses). In a strongly minimal set, the algebraic closure operator gives  a pregeometry, and there is an associated dimension assigned to any subset of $M$. The strongly minimal structure $\calM$ is {\em degenerate} if, for all $A\subset M$, we have $\acl(A)=\bigcup(\acl(a): a \in A)$. (Such strongly minimal sets are often called {\em trivial}, but we use {\em degenerate} to avoid confusion with the notion of a {\em trivial reduct}.) We say $\calM$ is {\em locally modular} if $\dim(A \cup B) +\dim(A\cap B)=\dim(A)+\dim(B)$ for all algebraically closed subsets $A,B$ of $M$ with $\dim(A\cap B)>0$. For more details, see Section 4.6 of \cite{hodges} or Appendix C of \cite{tz}.

We make the following  observations about reducts of strongly minimal sets. Recall that a graph $\Gamma$ is said to be {\em vertex-transitive} if for any vertices $x,y$ of $\Gamma$ there is $g\in \Aut(\Gamma)$ such that $x^g=y$. Also, the {\em degree} or {\em valency} of a vertex is its number of neighbours. A graph is said to be {\em regular} if all vertices have the same degree,  and to be {\em locally finite} if each vertex has finite degree.

\begin{remark} \rm \label{morenotes}
\begin{enumerate}
\item If $\calM$ is strongly minimal, then every  definable reduct $\calM'$ of $\calM$ is strongly minimal. If in addition $\calM$ is locally modular, then so is $\calM'$; see e.g. \cite[Corollary 6]{epp} or \cite[Proposition 6.3]{pillay}, which show  that any reduct of a superstable one-based theory of finite U-rank is one-based. Also, if $\calM$ is a locally modular but not degenerate strongly minimal set, then an infinite group is interpretable in $\calM$ -- this is well known and follows for example from Proposition 2.1 and Theorem 3(a) of \cite{hrush}. By for example Theorem 4.5 of \cite{zilber}, an infinite group can never be interpreted in a degenerate strongly minimal set. Thus, a definable reduct of a degenerate strongly minimal set is always degenerate.  
\item Any vertex-transitive graph of finite valency is a degenerate  strongly minimal structure; see Lemma 2.1 (ii), (iii) of \cite{blossier} -- the degenerateness follows from the quantifier-elimination in (ii).
\item The vertex-transitive graph $\calZ=({\mathbb Z},R)$, where $Rxy$ holds if and only if $|x-y|=1$, has infinitely many distinct definable reducts. Indeed, for any integer $n>1$, let $\calZ^{(n)}$ be the graph obtained from $\calZ$ by making two vertices adjacent if they are at distance $n$ in $\calZ$. 
Then $\calZ^{(n)}$ is the disjoint union of $n$ isomorphic copies of $\calZ$, so if $n\neq m$ then $\calZ^{(n)}$ and $\calZ^{(m)}$ are distinct definable reducts of $\calZ$. These are also distinct group-reducts.
\item Let $(T,R)$ be any regular tree of finite valency $t>2$. Let $T^{(2)}$ be the graph 
with vertex set $T$,
two vertices adjacent in $T^{(2)}$ if they are at distance 2 in $(T,R)$. Then $T^{(2)}$
 is a proper non-trivial definable reduct of $(T,R)$ (and also a proper non-trivial group reduct), and is the disjoint union of two graphs, each (denoted $\Gamma_{t-1,t}$ in Theorem~\ref{dist-trans} below) consisting of copies of 
$K_{t}$ joined in a treelike way, so that each vertex lies in $t$ copies of $K_{t}$. 

There is a further group-reduct of $(T,R)$ which does not correspond to a definable 
reduct: let $E$ be the equivalence relation on $T$, with two vertices equivalent if 
they are at even distance in $(T,R)$. There are exactly two $E$-classes, and 
$\Aut(T,E)\cong \Sym({\mathbb N}) \Wr C_2$. By strong minimality of $(T,R)$, the relation $E$ is not definable in  $(T,R)$.

We do not know whether $(T,R)$ has any other proper non-trivial definable reducts or group reducts (but see Theorem~\ref{dist-trans} below).
\item Let $F$ be a field, and $V$ be a vector space over $F$ with $|V|=\aleph_0$, viewed as a structure $\mathcal{V}$ in the language $(+,-,0,(f_a)_{a\in F})$, where $f_a$ is the unary function given by multiplication by $a$. Then $\mathcal{V}$ has quantifier
elimination, and is strongly minimal, and locally modular but not degenerate.
If $F$ is a finite field of prime order ${\mathbb F}_p$, then $\mathcal{V}$ is $\omega$-categorical so definable reducts and group reducts coincide. Its automorphism group ${\rm GL}(\aleph_0,p)$ has the affine group ${\rm AGL}(\aleph_0,p)=(V,+)\rtimes {\rm GL}(\aleph_0,p)$ as a closed supergroup, and the latter is maximal-closed in
$\Sym(V)$. Indeed, ${\rm AGL}(\aleph_0,p)$ is the automorphism group of an $\omega$-categorical strictly minimal set, and these are classified and shown to be the union of a sequence of `envelopes' 
 in \cite{chl} and \cite{zilber}. Thus, to see that ${\rm AGL}(\aleph_0,p)$ is maximal-closed in $\Sym(V)$ it is only necessary to observe that for $n>0$ the group ${\rm AGL}(n,p)$ is maximal in the finite symmetric group $S_{p^n}$ -- see the main theorem of \cite{mortimer} for this. In fact, this argument shows that if $q$ is a prime power then ${\rm AGL}(\aleph_0,q)$ has just finitely many closed supergroups. 

If $F$ is countably infinite and of characteristic 0 then $\mathcal{V}$ has a definable reduct with a single binary relation $R$, with $Rxy$ if and only if $x,y$ are nonzero and either $x+x=y$ or $y+y=x$. This graph is a disjoint union of infinitely many copies of the graph $\calZ$ of (2) above, along with an isolated vertex $\{0\}$, so has infinitely many definable (but isomorphic) reducts, and hence so does $\mathcal{V}$. These are also distinct group-reducts. In the case when $\dim V=1$, see also Proposition~\ref{valuation} for another construction of reducts. If $F$ is infinite of characteristic $p$ with prime subfield $\Ff_p$, then $\mathcal{V}$ has as a proper non-trivial definable reduct (and group-reduct) the structure of a vector space over $\Ff_p$, and intermediate fields give intermediate reducts.
\item Let $\calF=(F,+,\times)$ be an algebraically closed field. Then $\calF$ is strongly minimal and non locally modular. Of course, $\calF$ has proper group and definable  reducts of the form $(F,+)$ and $(F,\times)$.
If $\calF$ has characteristic zero, then we obtain infinitely many definable reducts as in (5): define $R$ on $F$, putting $Rxy$ if and only $x+x=y$ or $y+y=x$.

Suppose instead that ${\rm char}(\calF)=p>0$. Then there are definable and group reducts arising as above by viewing $F$ as an infinite-dimensional vector space (or affine space) over $\Ff_p$. To obtain further reducts, define a binary relation $R$ on $F$, putting $Rxy$ if and only if $x,y\not\in \{0,1\}$ and either $y=x^p$ or $x=y^p$. Then $(F,R)$ is a graph with finitely many cycles of each finite length (viewing connected components of size 1 and 2 as degenerate cycles), and infinitely many two-way infinite paths. Thus, $(F,R)$, and hence $(F,+,\times)$, has infinitely many distinct definable reducts and group-reducts, by (2).
\item As pointed out by Anand Pillay, another natural place to look for examples is the setting of $G$-sets -- see e.g. \cite[Exercise 4, p.169 and Exercise 4 p.177]{hodges}. Given a countably infinite group $G$ acting regularly (sharply 1-transitively) on a set $X$, we view $X$ as a structure with a unary function symbol for each group element, so $X$ is essentially a Cayley graph for $G$ with edges colored by group elements. This is a degenerate strongly minimal set (see the above Exercises in \cite{hodges}), and, as noted by Pillay,  proper definable reducts arise from proper non-trivial subgroups of $G$ (which will always exist). The automorphism group of the structure is exactly $G$ (in its regular action on $X$) and there will be a proper non-trivial group-reduct corresponding  to the action of $G\rtimes\Aut(G)$  on $X$. (Here, we identify  $X$ with $G$ by identifying some $x_0$ with $1$ and then each $x_0^g$ with $g$;  to define the action, if $x\in X, g\in G$ and $h\in \Aut(G)$ (so $(g,h)\in G\rtimes\Aut(G)$), we put $x^{(g,h)}=(x^h)g$, where $x^h$ denotes the image of $x$ -- viewed as an element of $G$ -- under the automorphism $h$.)
\end{enumerate}
\end{remark}

The two connected components of the graphs $T^{(2)}$ in Remark~\ref{morenotes} (4) above  belong to a larger  family. Let $k\geq 2$, $l\geq 2$, and let $T_{k+1,l}$ be the tree such that in the bipartition given by even distance, vertices in one part have valency $k+1$ and vertices in the other part have valency $l$. Form the graph $\Gamma_{k,l}$ whose vertex set is a part of the bipartition consisting of  $l$-valent elements of $T_{k+1,l}$, two adjacent if they are at distance two in $T_{k+1,l}$. The graphs $\Gamma_{k,l}$ are {\em distance-transitive}; that is, they are connected and if $u_1,u_2,v_1,v_2$ are vertices with $d(u_1,u_2)=d(v_1,v_2)$ then there is $g \in \Aut(\Gamma_{k,l})$ with $u_1^g=v_1$ and $u_2^g=v_2$. In fact these, and regular trees, are exactly the locally finite 
distance-transitive graphs classified in \cite{m}. Informally, $\Gamma_{k,l}$ consists of copies of the complete graph $K_{k+1}$ glued together in a treelike way, with $l$ copies containing each vertex (so the neighbourhood of a vertex consists of $l$ copies of $K_k$).

\begin{theorem} \label{dist-trans}
Let $k,l\in {\mathbb Z}$ with  $k\geq 2$ and $l\geq 3$. Then the graph $\Gamma_{k,l}$ has no proper non-trivial definable reducts.
\end{theorem}
\begin{proof}  Let $d$ denote the graph metric on $\Gamma_{k,l}$.
We first observe that by Remark~\ref{morenotes}(2), as $\Gamma_{k,l}$ is vertex-transitive of finite valency, it is strongly minimal and degenerate. It follows by Remark~\ref{morenotes}~(1) that any definable reduct $\calN$ of $\Gamma_{k,l}$ is also strongly minimal and degenerate. In particular, if $\calN$ is a non-trivial definable reduct then the algebraic closure in $\calN$ of each singleton $\{a\}$ has size greater than 1. (Indeed, otherwise, by transitivity and degenerateness every subset of $N$ (the domain of $\calN$) is algebraically closed in $\calN$, and it follows by strong minimality that $\Aut(\calN)=\Sym(N)$ and $\calN$ is a trivial reduct.) Thus, there is a formula $\phi(x,y)$ in the language of $\calN$ such that for each vertex $a$, the set $\phi(N,a)$ is finite of size greater than one. By distance-transitivity of $\Gamma_{k,l}$, there must be a finite subset
$\{n_1,\ldots,n_t\}$  of ${\mathbb N}$ such that $\phi(x,y)$ is equivalent to
$\bigvee_{i=1}^t d(x,y)=n_i$. We may suppose that $n_1<\ldots<n_t$, and put $n:=n_t$.

\medskip

{\em Claim.} If $u,v$ are vertices of $\Gamma_{k,l}$, then $d(u,v)\leq 2n$ if and only if 
$\calN\models \exists w \, (\phi(u,w)\wedge \phi(w,v))$.

\medskip

{\em Proof of Claim.}
The direction $\Leftarrow$ is immediate. For the direction $\Rightarrow$, let $m:=d(u,v)\leq 2n$. We may suppose $u\neq v$, as otherwise the result is immediate. If $m$ is even, say $m=2r$, let $z$ be the midpoint on the (unique) $uv$-geodesic, so $d(u,z)=d(v,z)=r$. Since $l\geq 3$ there is a neighbour $z'$ of $z$ at distance $r+1$ from $u$ and $v$ (so not adjacent to the neighbours of $z$ on the $zu$ or $zv$ geodesics). Hence we may choose $w$ so that $d(w,z)=d(w,z')+1=n-r$, so $d(w,u)=d(w,v)=n$, and $\calN\models \exists w(\phi(u,w)\wedge \phi(w,v))$.

Suppose instead that $m=2r+1>1$ is odd, and pick $z$ so that $d(z,u)=r$ and $d(z,v)=r+1$. Let $z'$ be the neighbour of $z$ on the $zv$-path, and $z''$ be a third vertex in the copy of $K_{k+1}$ containing $z$ and $z'$ (here we use that $k\geq 2$). Thus $d(u,z'')=d(v,z'')=r+1$. We now argue as in the last paragraph to find $w$ at distance $n-(r+1)$ from $z''$ and  at distance $n$ from $u$ and $v$. The argument if $m=1$ is similar. In that case, as $k\geq 2$ there is $z$ adjacent to $u$ and $v$ and we may choose $w$ distance $n-1$ from $z$ and  distance $n$ from $u$ and $v$. $\Box$

\medskip

If $v$ is a vertex of $\Gamma_{k,l}$ and $s$ is a positive integer, let $\gamma_s$ denote the number of vertices at distance $s$ from $v$ (which is clearly independent from the choice of $v$). Let $B_s(v)$ be the set of vertices at distance at most $s$ from $v$. Then as $\Gamma_{k,l}$ is locally finite, each $\gamma_s$ and $B_s(v)$ is finite; in fact,
$\gamma_s=lk^s(l-1)^{s-1}$ and 
$$|B_s(v)|=1+lk\frac{(l-1)^sk^s-1}{(l-1)k-1}.$$
There is also a fixed $\gamma$ such that if vertices $u$ and $v$ are adjacent in $\Gamma_{k,l}$ then 
$|B_{2n}(u)\cap B_{2n}(v)|= \gamma$ (and in fact, $\gamma=|B_{2n}(u)|-((l-1)k)^{2n}$).
To see that the adjacency relation of $\Gamma_{k,l}$ is $\emptyset$-definable in $N$, it suffices now to observe that if vertices $u,v$ are distinct and non-adjacent, then $|B_{2n}(u)\cap B_{2n}(v)|< \gamma$. This last assertion is easily seen. In fact, for such $u,v$, if $u'$ is chosen adjacent to $u$ and on the $uv$-geodesic then $B_{2n}(u)\cap B_{2n}(v)$ is a proper subset of $B_{2n}(u)\cap B_{2n}(u')$: the containment is clear by considering possible $uu'vx$-configurations, and if $d(u,x)=2n-1$ and $d(u',x)=2n$ then
$x\in (B_{2n}(u)\cap B_{2n}(u')) \setminus (B_{2n}(u)\cap B_{2n}(v))$.
\end{proof}

When the first draft of this paper was written, we did not know any example of a strongly minimal set which is  not $\omega$-categorical, has no non-trivial reducts (of either kind), and is not degenerate (but see the remarks about \cite{ks} below). 
An obvious potential example, now known by \cite{ks} to have no definable reducts,  is the 1-dimensional affine space ${\rm AG}_1({\mathbb Q})$, that is, the structure $({\mathbb Q},f)$ where $f$ is the ternary function defined by $f(x,y,z)=x-y+z$.
The automorphism group of $({\mathbb Q},f)$ is exactly the 1-dimensional affine group 
$\AGL(1,{\mathbb Q})=({\mathbb Q},+)\rtimes({\mathbb Q}^*,\cdot)$, and the following result shows that this group is {\em not} maximal-closed in $\Sym({\mathbb Q})$.

\begin{proposition}\label{valuation}
The structure $({\mathbb Q},f)$ has at least $\aleph_0$ proper non-trivial group-reducts which have pairwise incomparable automorphism groups and are not definable reducts.
\end{proposition}
\begin{proof}
For each prime $p$ let $v_p$ denote the $p$-adic valuation on the field $({\mathbb Q},+,\cdot)$. Define the ternary relation $C_p(x;y,z)$ on ${\mathbb Q}$, putting
$C_p(x;y,z)$ if and only if $v_p(x-y)<v_p(y-z)$. By \cite[Proposition 4.10]{ms}, 
$C_p$ is a $C$-relation on ${\mathbb Q}$ admitting $\AGL(1,{\mathbb Q})$ as a group of automorphisms. It is exactly the $C$-relation arising from the first description of $(M,D)$ in Example~\ref{mainex}, but working with sequences from $\{0,\ldots,p-1\}$ rather than just $\{0,1\}$. The group $\Aut({\mathbb Q},C_p)$ is a closed subgroup of $\Sym({\mathbb Q})$ containing $\AGL(1,{\mathbb Q})$. It is a Jordan group and clearly uncountable, so the containment is proper. 

We leave it to the reader the check that if $p,r$ are distinct primes then $\Aut({\mathbb Q},C_p)\neq \Aut({\mathbb Q},C_r)$. Finally, since $({\mathbb Q},f)$ is strongly minimal and so stable, and the structures  $({\mathbb Q},C_p)$ are unstable, the relations $C_p$ are not definable in $({\mathbb Q},f)$.
\end{proof}

As mentioned in Remark~\ref{notes}~(4), it is possible to define a $D$-relation $D_p$ analogously on $\PG(1,{\mathbb Q})$, using cross-ratio (see Theorem 30.4 in~\cite{an}). Thus, $\PGL(2,{\mathbb Q})$ is not maximal-closed, as 
$\PGL(2,{\mathbb Q})<\Aut({\mathbb Q},D_p)<\Sym({\mathbb Q})$. Again, distinct primes give distinct group-reducts. If $p=2$ then $({\mathbb Q},D_p)$
is exactly the structure $(M,D)$ of Theorem~\ref{reducts}, so the latter may be viewed as a group-reduct of the structure $\PG(1,{\mathbb Q})$ viewed in a language with a symbol $P_\lambda$ (for every $\lambda\in {\mathbb Q}\cup\{\infty\}$) holding of quadruples with cross-ratio $\lambda$. 

We conclude with some  questions.

\begin{question}\label{1,2}
\begin{enumerate}
\item Does the structure $({\mathbb Q},f)$ have any proper non-trivial definable reducts?
\item For $2\leq d\leq \aleph_0$ is the group $\AGL(d,{\mathbb Q})$  maximal-closed in
the symmetric group on the corresponding  affine space?  
\item Is there a non locally modular strongly minimal set with no proper non-trivial definable reduct (or group-reduct)? 
\item Do the graphs $\Gamma_{k,l}$ of Theorem~\ref{dist-trans} have any proper non-trivial group-reducts? Are their automorphism groups maximal in the symmetric group subject to being locally compact? (For a closed subgroup of $\Sym({\mathbb N})$, local compactness is equivalent to the pointwise stabiliser of some finite set having all orbits finite.)
\item Classify the definable reducts and group reducts of finite valency regular trees.
\end{enumerate}
\end{question}

Since the results and questions of this paper were communicated in a lecture, I. Kaplan and P. Simon \cite{ks} have answered (1) and (2). They use Jordan group methods to prove that $\AGL(n,{\mathbb Q})$ and $\PGL(n+1,{\mathbb Q})$ are maximal-closed for all $n$ with $2\leq n\leq \aleph_0$, and (with an observation of  Evans) deduce that $({\mathbb Q},f)$ has no proper non-trivial definable reducts. In particular, this yields countable maximal-closed groups. In an earlier version of this paper we also asked whether the groups $\PGL(n,{\mathbb Q})$ (for $3\leq n\leq \aleph_0$) are maximal-closed, and more generally whether $\Sym({\mathbb N})$ has any {\em countable} maximal-closed groups. However, both these questions are answered positively  in \cite{bogomolov}, where it is shown that if $H> \PGL(n,{\mathbb F})$ is closed ($F$ any field) and  contains a permutation taking a collinear triple to a non-collinear triple, then $H$ is the full symmetric group. 

We emphasise the broader question of finding more examples of non-oligomorphic maximal-closed subgroups of $\Sym({\mathbb N})$, and flag up the following question from \cite{mn}. Note that by \cite[Observation 3.3]{mn}, $\Sym({\mathbb N})$ is not the union of any chain of closed proper subgroups.

\begin{question}[Question 7.7 in~\cite{mn}] 
Is it true that every closed proper subgroup of $\Sym({\mathbb N})$ is contained in a 
maximal-closed subgroup of $\Sym({\mathbb N})$?
\end{question}

\begin{question} 
Show that there are $2^{\aleph_0}$ pairwise non-conjugate maximal-closed subgroups of $\Sym({\mathbb N})$.
\end{question}

The last question might be answered by the methods of 
\cite{bp1}. The automorphism groups of pairwise non-isomorphic `Henson digraphs' will be pairwise non-conjugate. It may be possible to show that they are maximal-closed using the approach developed in  \cite{bp1}. Indeed,  by \cite[Theorem 5.3]{nes}, for any Henson digraph determined by  a set of forbidden finite tournaments, the class of its finite subdigraphs, expanded in all possible ways by a total order, is a Ramsey class.

Observe that $2^{\aleph_0}$ is the maximal possible number of maximal-closed subgroups of $\Sym({\mathbb N})$, by easy counting. On the other hand,  $\Sym({\mathbb N})$ has $2^{2^{\aleph_0}}$ pairwise non-conjugate maximal subgroups, since stabilisers of ultrafilters on ${\mathbb N}$ are maximal and there are $2^{2^{\aleph_0}}$ distinct ultrafilters on ${\mathbb N}$ (see \cite[Theorem 6.4 and Corollary 6.5]{mn}).

\end{document}